\documentclass[11pt,reqno]{amsart}
\usepackage{amssymb,amscd,amsbsy}
\usepackage{amssymb,amscd,amsbsy,mathrsfs}
\newcommand{\lb}{\linebreak}

\renewcommand{\a}{\alpha}

\renewcommand{\d}{\delta}
\newcommand{\e}{\varepsilon}
\newcommand{\vk}{\varkappa}
\newcommand{\z}{\zeta}

\renewcommand{\l}{\lambda}

\newcommand{\s}{\sigma}
\renewcommand{\t}{\tau}

\newcommand{\f}{\varphi}
\renewcommand{\o}{\omega}

\newcommand{\D}{\Delta}
\renewcommand{\L}{\Lambda}

\renewcommand{\O}{\Omega}

\newcommand{\n}{\nabla}

\newcommand{\A}{{\mathcal A}}
\newcommand{\B}{{\mathcal B}}

\newcommand{\E}{{\mathcal E}}
\newcommand{\cd}{{\mathcal D}}
\newcommand{\F}{{\mathcal F}}
\newcommand{\Q}{{\mathcal Q}}
\newcommand{\h}{{\mathcal H}}
\newcommand{\I}{{\mathcal I}}

\newcommand{\K}{{\mathcal K}}

\newcommand{\cR}{{\mathcal R}}

\newcommand{\X}{{\mathcal X}}
\newcommand{\Y}{{\mathcal Y}}
\newcommand{\cU}{{\mathcal U}}
\newcommand{\V}{{\mathcal V}}

\newcommand{\C}{{\Bbb C}}
\newcommand{\T}{{\Bbb T}}
\newcommand{\pp}{{\Bbb P}}
\newcommand{\dd}{{\Bbb D}}
\newcommand{\R}{{\Bbb R}}
\newcommand{\Z}{{\Bbb Z}}

\newcommand{\0}{{\boldsymbol{0}}}

\newcommand{\bs}{\boldsymbol}

\newcommand{\m}{{\boldsymbol m}}
\newcommand{\bS}{{\boldsymbol S}}

\newcommand{\rf}[1]{(\ref{#1})}

\newcommand{\df}{\stackrel{\mathrm{def}}{=}}

\newcommand{\spn}{\operatorname{span}}
\newcommand{\supp}{\operatorname{supp}}
\newcommand{\clos}{\operatorname{clos}}

\newcommand{\const}{\operatorname{const}}

\newcommand{\eeq}{\end{equation}}
\newcommand{\beq}{\begin{equation}}
\newcommand{\bay}{\begin{eqnarray}}
\newcommand{\ba}{\begin{align*}}
\newcommand{\ea}{\end{align*}}
\newcommand{\ey}{\end{eqnarray}}
\newcommand{\bey}{\begin{eqnarray*}}
\newcommand{\eey}{\end{eqnarray*}}

\newcommand{\imp}{\Rightarrow}
\newcommand{\be}{\infty}

\newcommand{\bl}{\blacksquare}

\newcommand{\Pf}{{\bf Proof. }}

\newcommand{\ov}{\overline}

\newtheorem{thm}{\hspace{\parindent}Theorem}[section]

\newtheorem{cor}[thm]{\hspace{\parindent}Corollary}
\newtheorem{lem}[thm]{\hspace{\parindent}Lemma}

\theoremstyle{remark}

\newtheorem*{rem*}{Remark}

\newcommand\up{\upsilon}

\newcommand\fM{\frak M}
\newcommand\cZ{\mathcal{Z}}
\newcommand\dg{\frak D}

\begin{document}

\newcommand{\vse}{\vspace{.2in}}
\numberwithin{equation}{section}

\title{Operator H\"older--Zygmund functions}
\author{A.B. Aleksandrov and V.V. Peller}
\thanks{The first author is partially supported by RFBR grant 08-01-00358-a and by
Russian Federation presidential grant NSh-2409.2008.1;
the second author is partially supported by NSF grant DMS 0700995 and by ARC grant}

\newcommand{\mt}{{\mathcal T}}

\begin{abstract}
It is well known that a Lipschitz function on the real line does not have to be operator Lipschitz. We show that 
the situation changes dramatically if we pass to H\"older classes. Namely, we prove that
 if $f$ belongs to the H\"older class $\L_\a(\R)$ with $0<\a<1$, then $\|f(A)-f(B)\|\le\const\|A-B\|^\a$ for arbitrary self-adjoint operators $A$ and $B$. We prove a similar result for  functions $f$ in the Zygmund class $\L_1(\R)$: for arbitrary self-adjoint operators $A$ and $K$ we have $\|f(A-K)-2f(A)+f(A+K)\|\le\const\|K\|$. We also obtain analogs of this result for all H\"older--Zygmund classes $\L_\a(\R)$, $\a>0$. Then we find a sharp estimate for $\|f(A)-f(B)\|$ for functions $f$ of class $\L_\o\df\{f:~\o_f(\d)\le\const\o(\d)\}$ for an arbitrary modulus of continuity $\o$. In particular, we study moduly of continuity, for which $\|f(A)-f(B)\|\le\const\o(\|A-B\|)$ for self-adjoint $A$ and $B$, and for an arbitrary function $f$ in $\L_\o$. We obtain similar estimates for commutators $f(A)Q-Qf(A)$ and quasicommutators $f(A)Q-Qf(B)$. Finally, we estimate the norms of finite differences
$\sum\limits_{j=0}^m(-1)^{m-j}\left(\begin{matrix}m\\j\end{matrix}\right)f\big(A+jK\big)$ for $f$ in the class $\L_{\o,m}$ that is defined in terms of finite differences and a modulus continuity $\o$ of order $m$. We also obtaine similar results for unitary operators and for contractions.
\end{abstract}

\maketitle

\

\begin{center}
{\Large Contents}
\end{center}

\

\begin{enumerate}
\item[1.] Introduction \quad\dotfill \pageref{s1}
\item[2.]  Function spaces \quad\dotfill \pageref{s2}
\item[3.] Multiple operator integrals  \quad\dotfill \pageref{s3}
\item[4.] H\"older--Zygmund estimates for self-adjoint operators \quad\dotfill \pageref{s4}
\item[5.] The case of unitary operators \quad\dotfill \pageref{s5}
\item[6.]  The case of contractions \quad\dotfill \pageref{s6}
\item[7.]  Arbitrary moduli of continuity \quad\dotfill \pageref{s7}
\item[8.]  Operator continuous functions and operator moduli of continuity \quad\dotfill \pageref{s8}
\item[9.]  A universal family of self-adjoint operators\quad\dotfill \pageref{s9}
\item[10.]  Commutators and quasicommutators\quad\dotfill \pageref{s10}
\item[11.]  Higher order moduli of continuity\quad\dotfill \pageref{s11}
\item[] References \quad\dotfill \pageref{bibl}
\end{enumerate}

\

\setcounter{section}{0}
\section{\bf Introduction}
\setcounter{equation}{0}
\label{s1}

\

It is well known that a Lipschitz function on the real line is not necessarily {\it operator Lipschitz}, i.e.,
the condition 
$$
|f(x)-f(y)|\le\const|x-y|,\quad x,~y\in\R,
$$
does not imply that for self-adjoint operators
$A$ and $B$ on Hilbert space,
$$
\|f(A)-f(B)\|\le\const\|A-B\|.
$$
The existence of such functions was proved in \cite{F1}. Later in \cite{Pe1} necessary conditions were found for a function $f$ to be operator Lipschitz.
Those necessary conditions also imply that Lipschitz  functions do not have to be operator Lipschitz.
In particular, it was shown in \cite{Pe1} that an operator Lipschitz function must belong locally to the Besov space $B_1^1(\R)$ (see \S\,2 for an introduction to Besov spaces). Note that in \cite{Pe1} and \cite{Pe3} a stronger necessary condition was also obtained. 

It is also well known that a continuously differentiable function does not have to be operator differentiable. Moreover, the fact that $f$ is continuous differentiable does not imply that for bounded self-adjoint operators $A$ and $K$ the function
$$
t\mapsto f(A+tK)
$$
is differentiable. For $f$ to be operator differentiable it must satisfy the same necessary conditions \cite{Pe1}, \cite{Pe3}.
(Note that Widom posed in \cite{W} a problem entitled "When are differentiable functions differentiable?")

On the other hand it was proved in \cite{Pe1} and \cite{Pe3} that the condition that a function belongs to the Besov space $B_{\be1}^1(\R)$ is sufficient for operator Lipschitzness (as well as for operator differentiability). We also mention here the papers \cite{JW}, \cite{ABF}, \cite{KS1},
\cite{KS2}, \cite{KS3}, and \cite{KST} that study operator Lipschitz functions.

Many mathematicians working on such problems in perturbation theory believed that a similar situation occurs when 
considering H\"older classes of order $\a$ and operator H\"older classes of order $\a$, $0<\a<1$. In particular, Farforovskaya obtained in
\cite{F1}  the following estimate 
$$
\|f(A)-f(B)\|\le\const\|f\|_{\L_\a(\R)}\left(\log_2^2\frac{b-a}{\|A-B\|}+1\right)^\a\|A-B\|^\a
$$
for self-adjoint operators $A$ and $B$ with spectra in $[a,b]$ and for an arbitrary function $f$ in $\L_\a(\R)$, $0<\a<1$. She also obtained the same inequality for $\a=1$ and a Lipschitz function $f$ (see also \cite{F2}).

However, we show in this paper that {\it the situation changes  dramatically} if we consider H\"older classes $\L_\a(\R)$ with $0<\a<1$. 
In this case {\it H\"older functions are necessarily  operator H\"older}, i.e., the condition
\bay
\label{hol}
|f(x)-f(y)|\le\const|x-y|^\a,\quad x,~y\in\R,
\ey
implies that for self-adjoint operators $A$ and $B$ on Hilbert space,
\bay
\label{ohol}
\|f(A)-f(B)\|\le\const\|A-B\|^\a.
\ey
Note that the constant in \rf{ohol} depends not only on the constant in \rf{hol}, but also on $\a$ and must tend to infinity as the constant in
\rf{hol} is fixed and $\a$ goes to $1$.

We consider in this paper the same problem for the Zygmund class $\L_1(\R)$, i.e., the problem of whether a function $f$
in the Zygmund class $\L_1$ (i.e., $f$ is continuous and satisfies the inequality
$$
|f(x+t)-2f(x)+f(x-t)|\le\const|t|,\quad x,\,t\in\R)
$$
implies that $f$ is operator Zygmund, i.e.,
for arbitrary self-adjoint operators $A$ and $K$,
$$
\|f(A+K)-2f(A)+f(A-K)\|\le\const\|K\|.
$$
This problems was posed in \cite{F3}.

We show in this paper that the situation is the same as in the case of H\"older classes $\L_\a(\R)$, $0<\a<1$. Namely we prove that a Zygmund function must necessarily be operator Zygmund.

We also obtain similar results for the whole scale of H\"older--Zygmund classes $\L_\a(\R)$, $0<\a<\be$, of continuous functions $f$
satisfying
$$
\left|\sum_{k=0}^n(-1)^{n-k}\left(\begin{matrix}n\\k\end{matrix}\right)f(x+kt)\right|\le\const|t|^\a\quad(\mbox{here}\quad n-1\le\a <n).
$$
There are many natural equivalent (semi)norms on $\L_\a(\R)$,  for example, 
$$
\|f\|_{\L_\a(\R)}=\sup_{t\ne0}|t|^{-\a}\left|\sum_{k=0}^n(-1)^{n-k}\left(\begin{matrix}n\\k\end{matrix}\right)f(x+kt)\right|.
$$
The above results are obtained in \S\,4. In Sections 5 and 6 we obtain analogs of these results for unitary operators and for contractions.

In \S\,7 we estimate $\|f(A)-f(B)\|$ in terms of $\|A-B\|$ for functions $f$ of class $\L_\o$, (i.e., $|f(x)-f(y)|\le\const\o(|x-y|)\,$) for arbitrary moduli of continuity 
$\o$. In particular, we study those moduli of continuity, for which the fact that $f\in\L_\o$ implies that 
$$
\|f(A)-f(B)\|\le\const\o\big(\|A-B\|\big)
$$
for arbitrary self-adjoint operators $A$ and $B$. We compare this class of moduli of continuity with the class of moduli of continuity $\o$, for which
the Hilbert transform acts on $\L_\o$.

In \S\,8 we study the class of operator continuous functions and for a uniformly continuous function $f$ 
we introduce the operator modulus of continuity
$\O_f$. The material of \S\,9 is closely related to that of \S\,8. We construct a universal family $\{A_t\}_{t\ge0}$  of self-adjoint operators in the sense that
to compute $\O_f$ for arbitrary $f$, it suffices to consider the family $\{A_t\}_{t\ge0}$.

Section 10 is devoted to norm estimates for commutators $f(A)Q-Qf(A)$ and quasicommutators $f(A)Q-Af(B)$. We compare the operator modulus of continuity with several other moduli of continuity defined in terms of commutators and quasicommutators. 

In the last section we obtain norm estimates for finite differences 
\bay
\label{hod}
\sum_{j=0}^m(-1)^{m-j}\left(\begin{matrix}m\\j\end{matrix}\right)f\big(A+jK\big),
\ey
where $f$ belongs to the class $\L_{\o,m}$ that is defined in terms of finite differences and $\o$ is a modulus of continuity of order $m$.

In \S\,2 we collect necessary information on Besov classes (and in particular, the H\"older--Zygmund classes), and spaces $\L_\o$ and $\L_{\o,m}$.
In \S\,3 we give a brief introduction into double and multiple operator integrals.

Note that the main results of this paper were announced in \cite{AP1}. In \cite{AP2} we are going to study the problem of the behavior of functions of operators under perturbations of Schatten--von Neumann class $\bS_p$. We are going to study properties of functions of perturbed dissipative operators 
in \cite{AP4}, where we improve results of \cite{Nab}.

Finally, we would like mention that Farforovskaya and Nikolskaya have informed us recently that they had found another proof of the fact that
a H\"older function of order $\a$, $0<\a<1$, must be operator H\"older of order $\a$.

\

\section{\bf Function spaces}
\setcounter{equation}{0}
\label{s2}

\

{\bf 2.1. Besov classes.}
The purpose of this subsection is to give a brief introduction to the Besov spaces that play an important role in problems of perturbation theory.
We start with Besov spaces on the unit circle.

Let $1\le p,\,q\le\be$ and $s\in\R$. The Besov class $B^s_{pq}$ of functions (or
distributions) on $\T$ can be defined in the following way. Let $w$ be an infinitely differentiable function on $\R$ such
that
\bay
\label{w}
w\ge0,\quad\supp w\subset\left[\frac12,2\right],\quad\mbox{and} \quad w(x)=1-w\left(\frac x2\right)\quad\mbox{for}\quad x\in[1,2].
\ey

Consider the trigonometric polynomials $W_n$, and $W_n^\sharp$ defined by
$$
W_n(z)=\sum_{k\in\Z}w\left(\frac{k}{2^n}\right)z^k,\quad n\ge1,\quad W_0(z)=\bar z+1+z,\quad
\mbox{and}\quad W_n^\sharp(z)=\ov{W_n(z)},\quad n\ge0.
$$
Then for each distribution $f$ on $\T$,
$$
f=\sum_{n\ge0}f*W_n+\sum_{n\ge1}f*W^\sharp_n.
$$
The Besov class $B^s_{pq}$ consists of functions (in the case $s>0$) or distributions $f$ on $\T$
such that
\bay
\label{bes}
\big\{\|2^{ns}f*W_n\|_{L^p}\big\}_{n\ge1}\in\ell^q\quad\mbox{and}
\quad\big\{\|2^{ns}f*W^\sharp_n\|_{L^p}\big\}_{n\ge1}\in\ell^q.
\ey

To define a regularized de la Vall\'ee Poussin type kernel $V_n$, we define the $C^\be$ function $v$ on $\R$ by
\bay
\label{VP}
v(x)=1\quad\mbox{for}\quad x\in[-1,1]\quad\mbox{and}\quad v(x)=w(|x|)\quad\mbox{if}\quad |x|\ge1,
\ey
where $w$ is a function described in \rf{w}.
Then the trigonometric polynomial $V_n$ is defined by
$$
V_n(z)=\sum_{k\in\Z}v\left(\frac{k}{2^n}\right)z^k,\quad n\ge1.
$$

Besov classes admit many other descriptions. In particular, for $s>0$, the space $B^s_{pq}$ admits the
following characterization. A function $f\in L^p$ belongs to $B^s_{pq}$, $s>0$, if and only if
$$
\int_\T\frac{\|\D^n_\t f\|_{L^p}^q}{|1-\t|^{1+sq}}d\m(\t)<\be\quad\mbox{for}\quad q<\be
$$
and
\bay
\label{pbe}
\sup_{\t\ne1}\frac{\|\D^n_\t f\|_{L^p}}{|1-\t|^s}<\be\quad\mbox{for}\quad q=\be,
\ey
where $\m$ is normalized Lebesgue measure on $\T$, $n$ is an integer greater than $s$, and $\D_\t$, $\t\in\T$, is
the difference operator: 
$$
(\D_\t f)(\z)=f(\t\z)-f(\z), \quad\z\in\T.
$$

We use the notation $B_p^s$ for $B_{pp}^s$.

The spaces $\L_\a\df B_\be^\a$ form the {\it H\"older--Zygmund scale}. If $0<\a<1$, then $f\in\L_\a$ if and only if
$$
|f(\z)-f(\t)|\le\const|\z-\t|^\a,\quad\z,\,\t\in\T,
$$
while $f\in\L_1$ if and only if $f$ is continuous and
$$
|f(\z\t)-2f(\z)+f(\z\bar\t)|\le\const|1-\t|,\quad\z,\,\t\in\T.
$$
By \rf{pbe},  $\a>0$, $f\in\L_\a$ if and only if $f$ is continuous and
$$
|(\D^n_\t f)(\z)|\le\const|1-\t|^\a,
$$
where $n$ is a positive integer such that $n>\a$. 

Note that the (semi)norm of a function $f$ in $\L_\a$ is equivalent to
$$
\sup_{n\ge1}2^{n\a}\big(\|f*W_n\|_{L^\be}+\|f*W_n^\sharp\|_{L^\be}\big).
$$

We denote by $\l_\a$ the closure of the set of trigonometric polynomials in $\L_\a$. It is easy to see that $f$ belongs to $\l_\a$ if and only if
$$
\lim_{n\to\be}2^{n\a}\|f*W_n\|_{L^\be}=\lim_{n\to\be}2^{n\a}\|f*W^\sharp_n\|_{L^\be}=0.
$$
If $\a>0$, this is equivalent to the fact that
$$
\lim_{\t\to1}\frac{|(\D^n_\t f)(\z)|}{|1-\t|^\a}=0.
$$

It is well known that the dual space $(\l_\a)^*$ can be identified naturally with the Besov space $B_1^{-\a}$ with respect to the following pairing:
$$
\langle f,g\rangle=\sum_{n\in\Z}\hat f(n)\hat g(n)
$$
in the case when $g$ is trigonometric polynomial. It is also well known that the dual space $\big(B_1^{-\a}\big)^*$ can be identified naturally with 
the space $\L_\a$ with respect to the same pairing.

It is easy to see from the definition of Besov classes that the Riesz projection $\pp_+$,
$$
\pp_+f=\sum_{n\ge0}\hat f(n)z^n,
$$
is bounded on $B^s_{pq}$. Functions in $\big(B^s_{pq}\big)_+\df\pp_+B^s_{pq}$ admit a natural extension to analytic functions
in the unit disk $\dd$. It is well known that the functions in $\big(B^s_{pq}\big)_+$ admit the following description:
$$
f\in \big(B^s_{pq}\big)_+\Leftrightarrow
\int_0^1(1-r)^{q(n-s)-1}\|f^{(n)}_r\|^q_p\,dr<\be,\quad q<\be,
$$
and
$$
f\in \big(B^s_{p\be}\big)_+\Leftrightarrow
\sup_{0<r<1}(1-r)^{n-s}\|f^{(n)}_r\|_p<\be,
$$
where $f_r(\z)\df f(r\z)$ and $n$ is a nonnegative integer greater than $s$.

Let us proceed now to Besov spaces on the real line. We consider homogeneous Besov spaces 
$B_{pq}^s(\R)$ of functions (distributions) on $\R$.
We use the same function $w$
as in \rf{w} and define the functions $W_n$ and $W^\sharp_n$ on $\R$ by
$$
\F W_n(x)=w\left(\frac{x}{2^n}\right),\quad\F W^\sharp_n(x)=\F W_n(-x),\quad n\in\Z,
$$
where $\F$ is the {\it Fourier transform}:
$$
\big(\F f\big)(t)=\int_\R f(x)e^{-{\rm i}xt}\,dx,\quad f\in L^1.
$$

With every tempered distribution $f\in{\mathscr S}^\prime(\R)$ we 
associate a sequences $\{f_n\}_{n\in\Z}$,
$$
f_n\df f*W_n+f*W_n^\sharp.
$$
Initially we define the (homogeneous) Besov class $\dot B^s_{pq}(\R)$ as the set of all $f\in{\mathscr S}^\prime(\R)$
such that 
\bay
\label{Wn}
\{2^{ns}\|f_n\|_{L^p}\}_{n\in\Z}\in\ell^q(\Z).
\ey
According to this definition, the space $\dot B^s_{pq}(\R)$ contains all polynomials. Moreover, the distribution $f$ is defined by the sequence $\{f_n\}_{n\in\Z}$
uniquely up to a polynomial. It is easy to see that the series $\sum_{n\ge0}f_n$ converges in ${\mathscr S}^\prime(\R)$.
However, the series $\sum_{n<0}f_n$ can diverge in general. It is easy to prove that the
series $\sum_{n<0}f_n^{(r)}$ converges on uniformly $\R$ for each nonnegative integer
$r>s-1/p$. Note that in the case $q=1$ the series $\sum_{n<0}f_n^{(r)}$ converges uniformly, whenever $r\ge s-1/p$.

Now we can define the modified (homogeneous) Besov class $B^s_{pq}(\R)$. We say that a distribution $f$
belongs to $B^s_{pq}(\R)$ if $\{2^{ns}\|f_n\|_{L^p}\}_{n\in\Z}\in\ell^q(\Z)$ and
$f^{(r)}=\sum_{n\in\Z}f_n^{(r)}$ in the space ${\mathscr S}^\prime(\R)$, where $r$ is
the minimal nonnegative integer such that $r>s-1/p$ ($r\ge s-1/p$ if $q=1$). Now the function $f$ is determined uniquely by the sequence $\{f_n\}_{n\in\Z}$ up
to a polynomial of degree less that $r$, and a polynomial $\varphi$ belongs to $B^s_{pq}(\R)$
if and only if $\deg\varphi<r$.

We can also define de la Vall\'ee Poussin type functions $V_n$, $n\in\Z$, by
$$
\F V_n(x)=v\left(\frac{x}{2^n}\right),
$$
where $v$ is a function given by \rf{VP}.

We use the same notation $V_n$, $W_n$ and $W_n^\sharp$ for functions on $\T$ and on $\R$. This will not lead to a confusion.
For positive $n$ we can easily obtain the function $V_n$ on the circle from the corresponding function $V_n$ on the line. 
It suffices to consider the $2\pi$-periodic function
$$
\sum_{j\in\Z}V_n(x+2j\pi)
$$
and identify it with a function on $\T$. The same can be done with the functions $W_n$ and $W_n^\sharp$.

Besov spaces $B^s_{pq}(\R)$ admit equivalent definitions that are similar to those discussed above in the case of Besov
spaces of functions on $\T$. In particular, the H\"older--Zygmund classes $\L_\a(\R)\df B^\a_{\be}(\R)$, $\a>0$, can be described 
as the classes of continuous functions $f$ on $\R$ such that
$$
\big|(\D^m_tf)(x)\big|\le\const|t|^\a,\quad t\in\R,
$$
where the difference operator $\D_t$ is defined by
$$
(\D_tf)(x)=f(x+t)-f(x),\quad x\in\R,
$$
and $m$ is an integer greater than $\a$.

As in the case of functions on the unit circle, we can introduce the following equivalent (semi)norm on $\L_\a(\R)$:
$$
\sup_{n\in\Z}2^{n\a}\big(\|f*W_n\|_{L^\be}+\|f*W_n^\sharp\|_{L^\be}\big),\quad f\in\L_\a(\R).
$$

Consider now the class $\l_\a(\R)$, which is defined as the closure of the Schwartz class ${\mathscr S}(\R)$ in $\L_\a(\R)$.
The following result gives a description of $\l_\a(\R)$
for $\a>0$. We use the following notation: $C_0(\R)$ stands for the space of continuous functions $f$ on $\R$
such that $\lim\limits_{|x|\to\be}f(x)=0$; $f_n\df f*W_n+f*W_n^\sharp$.

\begin{thm}
\label{mH}
Let $\a>0$ and let $m$ be the integer such that $m-1\le\a<m$. Suppose that $f\in\L_\a(\R)$. The following are equivalent:

{\em(i)} $f\in\l_\a(\R)$;

{\em(ii)} $f_n\in C_0(\R)$ for every $n\in\Z$ and 
$$
\lim_{|n|\to\be}2^{n\a}\|f_n\|_{L^\be}=0;
$$

{\em(iii)} the following equalitites hold:
$$
\lim_{t\to0}|t|^{-\a}\big(\D_t^mf\big)(x)=0\quad\mbox{uniformly in}\quad x\in\R,
$$
$$
\lim_{|t|\to\be}|t|^{-\a}\big(\D_t^mf\big)(x)=0\quad\mbox{uniformly in}\quad x\in\R,
$$
and
$$
\lim_{|x|\to\be}|t|^{-\a}\big(\D_t^mf\big)(x)=0\quad\mbox{uniformly in}\quad t\in\R\setminus\{0\}.
$$
\end{thm}

\Pf (ii)$\imp$(i). It follows from the definition of $\L_\a(\R)$ in terms of convolutions with $W_n$ and $W_n^\sharp$ that 
$$
\lim\left\|f-\sum_{n=-N}^Nf_n\right\|_{\L_\a(\R)}=0.
$$
Thus it suffices to prove that $f_n\in\l_\a(\R)$. However, this is a consequence of the following easily verifiable fact:
$$
\lim_{\e\to0}\,\sup_{x\in\R}\,\left|\left(e^{-\e^2x^2}f_n(x)\right)^{(j)}-f_n^{(j)}(x)\right|=0\quad\mbox{for all}\quad j\ge0.
$$

The implication (i)$\imp$(iii) follows very easily from the fact that (iii) holds for all functions in ${\mathscr S}(\R)$ which can easily be established.

It remains to show that (iii) implies (ii). Consider the function $Q_n$ defined by 
\bay
\label{Qn}
Q_n(t)=\sum_{k=1}^m(-1)^{k-1}\left(\begin{matrix}m\\k\end{matrix}\right)\frac1kV_n\left(\frac tk\right).
\ey
It is easy to see that
\begin{align}
\label{Qn0}
f(x)-\big(f*Q_n\big)(x)&=f(x)-\int_\R f(x-t)\sum_{k=1}^m(-1)^{k-1}\left(\begin{matrix}m\\k\end{matrix}\right)\frac1kV_n\left(\frac {t}k\right)\,dt\nonumber\\[.2cm]
&=f(x)+\int_\R \sum_{k=1}^m(-1)^{k}\left(\begin{matrix}m\\k\end{matrix}\right)f(x-kt)V_n(t)\,dt\nonumber\\[.2cm]
&=\int_\R\big(\D_{-t}^mf\big)(x)V_n(t)\,dt.
\end{align}
Hence,
\begin{align*}
2^{\a n}\big\|f-f*Q_n\big\|_{L^\be}&=\sup_{x\in\R}
2^{\a n}\left|\int_\R\big(\D^m_{-t}f\big)(x)V_n(t)\,dt\right|\\[.2cm]
&=\sup_{x\in\R}\left|\int_\R\frac{\big(\D^m_{-2^{-n}t}f\big)(x)}{|t|^\a2^{-\a n}}V(t)|t|^\a\,dt\right|\to0\quad\mbox{as}\quad|n|\to\be
\nonumber
\end{align*}
by the Lebesgue dominant convergence theorem.

Let us observe now that $\supp\F Q_n\subset\big[-2^{n+1},2^{n+1}\big]$, and so
\begin{align*}
\|f-f*V_n\|_{L^\be}&=\|f-f*Q_{n-1}-(f-f*Q_{n-1})*V_n\|_{L^\be}\\[.2cm]
&\le\|f-f*Q_{n-1}\|_{L^\be}+\|(f-f*Q_{n-1})*V_n\|_{L^\be}\\[.2cm]
&\le\const\|f-f*Q_{n-1}\|_{L^\be}
\end{align*}
which immediately implies that
$$
\lim_{|n|\to\be}2^{\a n}\|f_n\|_{L^\be}=0.
$$
Similarly, we can prove that $f-f*Q_n\in C_0(\R)$ and $f_n\in C_0(\R)$. $\bl$

The dual space $\big(\l_\a(\R)\big)^*$ to $\l_\a(\R)$ can be identified in a natural way with $B_1^{-\a}(\R)$ with respect to the pairing
$$
\langle f,g\rangle\df \lim_{N\to\be}\sum_{n=-N}^N\int_\R\big( \F(f_n)\big)(t)\big(\F g\big)(t)\,dt,\quad f\in\l_\a(\R),~g\in B_1^{-\a}(\R).
$$
The dual space $\big(B_1^{-\a}(\R)\big)^*$ to $B_1^{-\a}(\R)$ can be identified with $\L_\a(\R)$ with respect to the same pairing.

We refer the reader to \cite{Pee} and \cite{Pe4} for more detailed information on Besov spaces. 

We conclude this subsection with the following result that will be used in \S\,4.

\begin{thm}
\label{prodol}
Let $\a>0$. Then
for each $\varepsilon>0$ and each function
$f\in\Lambda_\alpha(\R)$
there exists a function $g\in\Lambda_\alpha(\R)$ with compact support such that
$f(t)=g(t)$ for $t\in[0,1]$ and
$$
\|g\|_{\Lambda_\alpha}\le\const\|f\|_{\Lambda_\alpha}+\varepsilon,
$$
where the constant can depend only on $\a$.
\end{thm}

To prove Theorem \ref{prodol}, we use the well-known fact that if $\f$ and $f$ are functions in $\L_\a(\R)$ and $\f$ has compact support, then $\f f\in \L_\a(\R)$.
We refer the reader to \cite{T}, Section 4.5.2 for the proof.

\begin{lem}
\label{pol}
Let $\a>0$ and let $P$ be a polynomial whose degree is at most $\a$. Then for an arbitrary $\e>0$ there exists
a function $f\in\L_\a(\R)$ with compact support such that
$$
f\big|[0,1]=P\big|[0,1]\quad\mbox{and}\quad\|f\|_{\L_\a(\R)}<\e.
$$
\end{lem}

\Pf It suffices to consider the case when $P(x)=x^n$ with $n\le\a$. Assume first that $n<\a$.
Let $g$ be an arbitrary function in $\L_\a(\R)$ with compact support and such that $g(x)=x^n$ for $x\in[0,1]$.
For $t\in(0,1)$, we define the function $g_t$ by
$$
g_t(x)=t^{-n}g(tx).
$$
It is easy to see that $g_t(x)=x^n$ for $x\in[0,1]$ and
$$
\|g_t\|_{\Lambda_\alpha(\R)}=t^{\alpha-n}\|g\|_{\Lambda_\alpha(\R)}\to0\quad\mbox{as}\quad t\to0.
$$

Suppose now that $\a$ is an integer and $n=\a$. It is well known that the function $h$ defined by $h(x)=x^n\log|x|$ belongs to
$\L_n(\R)$. Multiplying it by a suitable function in $\L_n(\R)$ with compact support, we obtain a function $g\in\L_n(\R)$ with compact support such that
$g(x)=x^n\log|x|$ for $x\in[0,1]$. For $t\in(0,1)$, we define the function $g_t$ by
$$
g_t(x)=(t^{-n}g(tx)-g(x))/\log t.
$$
Then $g_t(x)=x^n$ for $x\in[0,1]$ and
$$
\|g_t\|_{\Lambda_n(\R)}\le2|\log t|^{-1}\|g\|_{\Lambda_n(\R)}\to0\quad\mbox{as}\quad t\to0.\quad\bl
$$

\medskip

{\bf Proof of Theorem \ref{prodol}.}  Let $\f$ be a function in $\L_\a(\R)$ with compact support.
We fix a subset $\D$ of $[0,1]$ that has $n$ elements, where
$n$ is the largest integer such that $n\le\a+1$.
It follows from the closed graph theorem that
$\|\varphi f\|_{\Lambda_\alpha}\le C(\varphi,\alpha,\D)\|f\|_{\Lambda_\alpha}$
for every $f\in\Lambda_\alpha$ that vanishes on $\D$.
It remains to observe that an arbitrary function in $\Lambda_\alpha$ can be represented as 
the sum of a polynomial of degree at most  $\alpha$ and a function $\Lambda_\alpha$
vanishing on $\D$. $\bl$

\medskip

{\bf 2.2. Spaces $\bs{\L_\o}$.} Let $\o$ be a modulus of continuity, i.e., $\o$ is a nondecreasing continuous function on $[0,\be)$
such that $\o(0)=0$, $\o(x)>0$ for $x>0$, and
$$
\o(x+y)\le\o(x)+\o(y),\quad x,~y\in[0,\be).
$$
We denote by $\L_\o(\R)$ the space of functions on $\R$ such that
$$
\|f\|_{\L_\o(\R)}\df\sup_{x\ne y}\frac{|f(x)-f(y)|}{\o(|x-y|)}.
$$
We also consider in this paper the spaces  $\L_\o$ of functions on the unit circle and $\big(\L_\o\big)_+$ of functions analytic in the unit disc
that can be defined in a similar way.

\begin{thm}
\label{Vn}
There exists a constant $c>0$ such that for an arbitrary
modulus of continuity $\o$ and for an arbitrary function $f$ in $\L_\o(\R)$, 
the following inequality holds:
\bay
\label{VPn}
\|f-f*V_n\|_{L^\be}\le c\,\o\big(2^{-n}\big)\|f\|_{\L_\o(\R)},\quad n\in\Z.
\ey
\end{thm}

\Pf We have
\begin{align*}
\big|f(x)-\big(f*V_n\big)(x)\big|&=2^{n}\left|\int_\R\big(f(x)-f(x-y)\big)V\left(2^ny\right)\,dy\right|\\[.2cm]
&\le2^{n}\|f\|_{\L_\o(\R)}\int_\R\o(|y|)\,\left|V\left(2^ny\right)\right|\,dy\\[.2cm]
&=2^{n}\|f\|_{\L_\o(\R)}\int_{-2^{-n}}^{2^{-n}}\o(|y|)\,\left|V\left(2^ny\right)\right|\,dy\\[.2cm]
&+2^{n+1}\|f\|_{\L_\o(\R)}\int_{2^{-n}}^\be\o(y)\,\left|V\left(2^ny\right)\right|\,dy.
\end{align*}
Clearly,
$$
2^{n}\int_{-2^{-n}}^{2^{-n}}\o(|y|)\,\left|V\left(2^ny\right)\right|\,dy
\le\o\big(2^{-n}\big)\|V\|_{L^1}.
$$
On the other hand, keeping in mind the obvious inequality 
$2^{-n}\o(y)\le2y\o\big(2^{-n}\big)$ for $y\ge2^{-n}$, we obtain
\begin{align*}
2^{n+1}\int_{2^{-n}}^\be\o(y)\,\left|V\left(2^ny\right)\right|\,dy&\le
4\cdot2^{2n}\o\big(2^{-n}\big)\int_{2^{-n}}^\be y\,\left|V\left(2^ny\right)\right|\,dy\\[.2cm]
&=4\,\o\big(2^{-n}\big)\int_{1}^\be y\,\left|V\left(y\right)\right|\,dy
\le\const\o\big(2^{-n}\big).
\end{align*}
This proves \rf{VPn}. $\bl$

\medskip

{\bf Remark.} A similar inequality holds for functions $f$ on $\T$ of class $\L_\o$:
$$
\|f-f*V_n\|_{L^\be}\le \const\,\o\big(2^{-n}\big)\|f\|_{\L_\o},\quad n>0.
$$
To prove it, it suffices to identify $f$ with a $2\pi$-periodic function on $\R$ and apply Theorem \ref{Vn}.

\begin{cor}
\label{Wnn}
Let $f\in\L_\o(\R)$. Then
$$
\|f*W_n\|_{L^\be}\le\const\o\big(2^{-n}\big)\|f\|_{\L_\o(\R)},\quad n\in\Z,
$$
and
$$
\|f*W^\sharp_n\|_{L^\be}\le\const\o\big(2^{-n}\big)\|f\|_{\L_\o(\R)},\quad n\in\Z.
$$
\end{cor}

\medskip

{\bf 2.3. Spaces $\bs{\L_{\o,m}}$.}
We proceed now to moduli of continuity of higher order. For a continuous function $f$ on $\R$, we define the $m$th modulus of continuity $\o_{f,m}$ of $f$ by 
$$
\o_{f,m}(x)=\sup_{\{h:0\le h\le x\}}\big\|\D_h^mf\big\|_{L^\be}=\sup_{\{h:0\le|h|\le x\}}\big\|\D_h^mf\big\|_{L^\be},\quad x>0.
$$

The following elementary formula  can easily be verified by induction:
\bay
\label{ind}
\big(\D_{2h}^mf\big)(x)=\sum\limits_{j=0}^m{m\choose j}\big(\D_{h}^mf\big)(x+jh).
\ey
It follows from \rf{ind} that 
$\o_{f,m}(2x)\le 2^m\o_{f,m}(x),\quad x>0$.

Suppose now that 
$\o$ is a nondecreasing function on $(0,\be)$ such that 
\bay
\label{on}
\lim_{x\to0}\o(x)=0\quad\mbox{and}\quad
\o(2x)\le2^m\o(x)\quad\mbox{for}\quad x>0.
\ey
 It is easy to see that in this case 
 \bay
 \label{udv}
 \o(tx)\le2^mt^m\o(x),\quad\mbox{for all}\quad x>0\quad\mbox{and}\quad t>1.
 \ey
Denote by $\L_{\o,m}(\R)$ the set of continuous functions $f$ on $\R$
satisfying
$$
\|f\|_{\L_{\o,m}(\R)}\df\sup\limits_{t>0}\frac{\|\D^m_tf\|_{L^\infty}}{\o(t)}<+\infty.
$$

\begin{thm}
\label{mnn}
There exists $c>0$ such that for an arbitrary
nondecreasing function $\o$ on $(0,\be)$ satisfying {\em\rf{on}} and
for an arbitrary function $f\in\L_{\o,m}(\R)$, the following inequality holds:
$$
\|f-f*V_n\|_{L^\be}\le c\,\o\big(2^{-n}\big)\|f\|_{\L_{\o,m}(\R)},\quad n\in\Z.
$$
\end{thm}




\Pf Consider the function $Q_n$ defined by \rf{Qn}. Applying formula \rf{Qn0}, we obtain
\begin{align*}
\big|f(x)-\big(f*Q_n\big)(x)\big|&=\left|\int_\R\big(\D_{-t}^mf\big)(x)V_n(t)\,dt\right|\le\|f\|_{\L_{\o,m}(\R)}\int_\R\o(|t|)|V_n(t)|\,dt.
\end{align*}
It follows from \rf{udv} that
\begin{align*}
\int_\R\o(|t|)|V_n(t)|\,dt&=\int_{-2^{-n}}^{2^n}(|t|)|V_n(t)|\,dt+2^{n+1}\int_{2^{-n}}^\be\o(t)|V\big(2^nt\big)|\,dt\\[.2cm]
&\le\|V_n\|_{L^1}\,\o\big(2^{-n}\big)+2^{n+1}\cdot2^{m(n+1)}\o\big(2^{-n}\big)
\int_{2^{-n}}^\be t^m|V\big(2^nt\big)|\,dt\\[.2cm]
&=\|V\|_{L^1}\,\o\big(2^{-n}\big)+2^{m+1}\,\o\big(2^{-n}\big)\int_{1}^\be t^m|V(t)|\,dt
\le\const\o\big(2^{-n}\big).
\end{align*}Summarizing the above estimates, we obtain
$$
\|f-f*Q_n\|_{L^\be}\le\const\,\o\big(2^{-n}\big)\|f\|_{\L_{\o,m}(\R)}.
$$

As in the proof of Theorem \ref{mH}, we have
\begin{align*}
\|f-f*V_n\|_{L^\be}&=\|f-f*Q_{n-1}-(f-f*Q_{n-1})*V_n\|_{L^\be}\\[.2cm]
&\le\|f-f*Q_{n-1}\|_{L^\be}+\|(f-f*Q_{n-1})*V_n\|_{L^\be}\\[.2cm]
&\le\const\|f-f*Q_{n-1}\|_{L^\be}\le\const\,\o\big(2^{-n}\big)\|f\|_{\L_{\o,m}(\R)}.\quad\bl
\end{align*}

\begin{cor}
\label{Wnm}
Let $f\in\L_{\o,m}(\R)$. Then
$$
\|f*W_n\|_{L^\be}\le\const\o\big(2^{-n}\big)\|f\|_{\L_\o(\R)},\quad n\in\Z,
$$
and
$$
\|f*W^\sharp_n\|_{L^\be}\le\const\o\big(2^{-n}\big)\|f\|_{\L_\o(\R)},\quad n\in\Z.
$$
\end{cor}

{\bf Remark.} As in the case $m=1$, a similar result holds for the space $\L_{\o,m}$ of functions on the unit circle, which consists of 
continuous $f$ functions such that
$$
\|f\|_{\L_{\o,m}}\df\sup_{\t\ne1}\frac{\big|\big(\D_\t^mf\big)(\z)\big|}{\o(|1-\t|)}<\be.
$$
Again, identifying a function $f$ in $\L_{\o,m}$ with a $2\pi$-periodic function on $\R$, we can see that
$$
\|f-f*V_n\|_{L^\be}\le\const\,\o\big(2^{-n}\big)\|f\|_{\L_{\o,m}},\quad n>0.
$$

\

\section{\bf Multiple operator integrals}
\setcounter{equation}{0}
\label{s3}

\

{\bf 3.1. Double operator integrals.}
In this section we give a brief introduction in the theory of double operator integrals
developed by Birman and Solomyak in \cite{BS1}, \cite{BS2}, and \cite{BS3}, see also their survey \cite{BS5}.

Let $(\X,E_1)$ and $(\Y,E_2)$ be spaces with spectral measures $E_1$ and $E_2$
on Hilbert spaces $\h_1$ and $\h_2$. Let us first define double operator integrals
\bay
\label{doi}
\int\limits_{\X}\int\limits_{\Y}\Phi(x,y)\,d E_1(x)\,Q\,dE_2(y),
\ey
for bounded measurable functions $\Phi$ and operators $Q:\h_2\to\h_1$
of Hilbert--Schmidt class $\bS_2$. Consider the set function $F$ whose values are orthogonal
projections on the Hilbert space $\bS_2(\h_2,\h_1)$ of Hilbert--Schmidt operators from $\h_2$ to $\h_1$, which is defined 
on measurable rectangles by
$$
F(\D_1\times\D_2)Q=E_1(\D_1)QE_2(\D_2),\quad Q\in\bS_2(\h_2,\h_1),
$$ 
$\D_1$ and $\D_2$ being measurable subsets of $\X$ and $\Y$. Note that left multiplication by $E_1(\D_1)$
obviously commutes with right multiplication by $E_2(\D_2)$.

 It was shown in \cite{BS4} that $F$ extends to a spectral measure on 
$\X\times\Y$. If $\Phi$ is a bounded measurable function on $\X\times\Y$, we define
$$
\int\limits_{\X}\int\limits_{\Y}\Phi(x,y)\,d E_1(x)\,Q\,dE_2(y)=
\left(\,\,\int\limits_{\X_1\times\X_2}\Phi\,dF\right)Q.
$$
Clearly,
$$
\left\|\,\,\int\limits_{\X}\int\limits_{\Y}\Phi(x,y)\,dE_1(x)\,Q\,dE_2(y)\right\|_{\bS_2}
\le\|\Phi\|_{L^\be}\|Q\|_{\bS_2}.
$$


If the transformer
$$
Q\mapsto\int\limits_{\X}\int\limits_{\Y}\Phi(x,y)\,d E_1(x)\,Q\,dE_2(y)
$$
maps the trace class $\bS_1$ into itself, we say that $\Phi$ is a {\it Schur multiplier of $\bS_1$ associated with 
the spectral measures $E_1$ and $E_2$}. In
this case the transformer
\bay
\label{tra}
Q\mapsto\int\limits_{\Y}\int\limits_{\X}\Phi(x,y)\,d E_2(y)\,Q\,dE_1(x),\quad Q\in \bS_2(\h_1,\h_2),
\ey
extends by duality to a bounded linear transformer on the space of bounded linear operators from $\h_1$ to $\h_2$
and we say that the function $\Psi$ on $\X_2\times\X_1$ defined by 
$$
\Psi(y,x)=\Phi(x,y)
$$
is {\it a Schur multiplier of the space of bounded linear operators associated with $E_2$ and $E_1$}.
We denote the space of such Schur multipliers by $\fM(E_2,E_1)$.


To state a very important formula by Birman and Solomyak, we consider for a continuously differential function
$f$ on $\R$, the divided difference $\dg f$,
$$
(\dg f)(x,y)\df\frac{f(x)-f(y)}{x-y},\quad x\ne y,\quad (\dg f)(x,x)\df f'(x)\quad x,\,y\in\R.
$$
Birman in Solomyak proved in \cite{BS3} that if
 $A$ is a self-adjoint operator (not necessarily bounded),
$K$ is a bounded self-adjoint operator, and
$f$ is a continuously differentiable 
function on $\R$ such that
$\dg f\in\fM(E_{A+K},E_A)$, then
\bay
\label{BSF}
f(A+K)-f(A)=\iint\limits_{\R\times\R}\big(\dg f\big)(x,y)\,dE_{A+K}(x)K\,dE_A(y)
\ey
and
$$
\|f(A+K)-f(A)\|\le\const\|\dg f\|_{\fM}\|K\|,
$$
where $\|\dg f\|_{\fM}$ is the norm of $\dg f$ in $\fM(E_{A+K},E_A)$.




A similar formula and
similar results also hold for unitary operators, in which case we have to integrate the divided
difference $\dg f$ of a function $f$ on the unit circle with respect to the spectral measures of the corresponding operator integrals.

It is easy to see that if a function $\Phi$ on $\X\times\Y$ belongs to the {\it projective tensor
product}
$L^\be(E_1)\hat\otimes L^\be(E_2)$ of $L^\be(E_1)$ and $L^\be(E_2)$ (i.e., $\Phi$ admits a representation
\bay
\label{ptp}
\Phi(x,y)=\sum_{n\ge0}\f_n(x)\psi_n(y),
\ey
where $\f_n\in L^\be(E_1)$, $\psi_n\in L^\be(E_2)$, and
\bay
\label{ptpn}
\sum_{n\ge0}\|\f_n\|_{L^\be}\|\psi_n\|_{L^\be}<\be),
\ey
then $\Phi\in\fM(E_1,E_2)$, i.e., $\Phi$ is a Schur multiplier of the space of bounded linear operators. For such functions $\Phi$ we have
$$
\int\limits_\X\int\limits_\Y\Phi(x,y)\,d E_1(x)Q\,dE_2(y)=
\sum_{n\ge0}\left(\,\int\limits_\X \f_n\,dE_1\right)Q\left(\,\int\limits_\Y \psi_n\,dE_2\right).
$$ 
Note that if $\Phi$ belongs to the projective tensor product $L^\be(E_1)\hat\otimes L^\be(E_2)$, its norm in $L^\be(E_1)\hat\otimes L^\be(E_2)$
is, by definition, the infimum of the  left-hand side of \rf{ptpn} over all representations \rf{ptp}.


More generally, $\Phi$ is a Schur multiplier  if $\Phi$ 
belongs to the {\it integral projective tensor product} $L^\be(E_1)\hat\otimes_{\rm i}L^\be(E_2)$ of $L^\be(E_1)$ and $L^\be(E_2)$, i.e., $\Phi$ admits a representation
$$
\Phi(x,y)=\int_\O \f(x,\o)\psi(y,\o)\,d\s(\o),
$$
where $(\O,\s)$ is a measure space, $\f$ is a measurable function on $\X\times \O$,
$\psi$ is a measurable function on $\Y\times \O$, and
$$
\int_\O\|\f(\cdot,\o)\|_{L^\be(E_1)}\|\psi(\cdot,\o)\|_{L^\be(E_2)}\,d\s(\o)<\be.
$$
If $\Phi\in L^\be(E_1)\hat\otimes_{\rm i}L^\be(E_2)$, then
\bay
\label{iptr}
\iint\limits_{\X\times\Y}\!\Phi(x,y)\,d E_1(x)\,Q\,dE_2(y)\!=\!\!
\int\limits_\O\!\left(\,\int\limits_\X \f(x,\o)\,dE_1(x)\!\right)\!Q\!
\left(\,\int\limits_\Y \psi(y,\o)\,dE_2(y)\!\right)\!d\s(\o).
\ey
Clearly, the function 
$$
\o\mapsto \left(\,\int\limits_\X \f(x,\o)\,dE_1(x)\right)Q
\left(\,\int\limits_\Y \psi(y,\o)\,dE_2(y)\right)
$$
is weakly measurable and
$$
\int\limits_\O\left\|\left(\,\int\limits_\X \f(x,\o)\,dE_1(x)\right)Q
\left(\,\int\limits_\Y \psi(y,\o)\,dE_2(y)\right)\right\|\,d\s(\o)<\be.
$$

It turns out that all Schur multipliers of the space of bounded linear operators can be obtained in this way (see \cite{Pe1}).

In connection with the Birman--Solomyak formula, it is important to obtain sharp estimates of divided differences in integral projective tensor products of $L^\be$ spaces. It was shown in \cite{Pe1} that if $f$ is a trigonometric polynomial of degree $d$, then 
\bay
\label{Bp}
\big\|\dg f\big\|_{C(\T)\hat\otimes C(\T)}\le\const d\,\|f\|_{L^\be}.
\ey
On the other hand, it was shown in \cite{Pe3} that if $f$ is a bounded function on $\R$ whose Fourier transform is supported on $[-\s,\s]$
(in other words, $f$ is an entire function of exponential type at most $\s$ that is bounded on $\R$), then $\dg f\in L^\be\hat\otimes_{\rm i}L^\be$
and
\bay
\label{Be}
\big\|\dg f\big\|_{L^\be\hat\otimes_{\rm i} L^\be}\le\const \s\|f\|_{L^\be(\R)}.
\ey
Note that inequalities \rf{Bp} and \rf{Be} were proved in \cite{Pe1} and \cite{Pe3} under the assumption that the Fourier transform of $f$ is supported
on $\Z_+$ (or $\R_+$); however it is very easy to deduce the general results from those partial cases.

Inequalities \rf{Bp} and \rf{Be} led in \cite{Pe1} and \cite{Pe3} to the fact that functions in $B_{\be1}^1$ and $B_{\be1}^1(\R)$ are operator Lipschitz.

It was observed in \cite{Pe3} that it follows from \rf{BSF} and \rf{Be} that if $f$ is an entire function of exponential type at most $\s$ that is bounded on $\R$, and $A$ and $B$ are self-adjoint operators with bounded $A-B$, then
$$
\|f(A)-f(B)\|\le\const\s\|f\|_{L^\be}\|A-B\|.
$$
Actually, it turns out that the last inequality holds with constant equal to 1. This will be proved in \cite{AP3}.

\medskip

{\bf 3.2. Multiple operator integrals.}  The approach by Birman and Solomyak to double operator integrals does not generalize to the case of
multiple operator integrals. However, formula \rf{iptr} suggests an approach to multiple operator integrals that is based on integral projective tensor products. This approach was given in \cite{Pe5}.

To simplify the notation, we consider here the case of triple operator integrals; the case of arbitrary multiple operator integrals can be treated in the same way.

Let $(\X,E_1)$, $(\Y,E_2)$, and $(\cZ,E_3)$
be spaces with spectral measures $E_1$, $E_2$, and $E_3$ on Hilbert spaces $\h_1$, $\h_2$, and $\h_3$. Suppose that
$\Phi$ belongs to the integral projective tensor product
$L^\be(E_1)\hat\otimes_{\rm i}L^\be(E_2)\hat\otimes_{\rm i}L^\be(E_3)$, i.e., $\Phi$ admits a representation
\bay
\label{ttp}
\Phi(x,y,z)=\int_\O \f(x,\o)\psi(y,\o)\chi(z,\o)\,d\s(\o),
\ey
where $(\O,\s)$ is a measure space, $\f$ is a measurable function on $\X\times \O$,
$\psi$ is a measurable function on $\Y\times \O$, $\chi$ is a measurable function on $\cZ\times \O$,
and
$$
\int_\O\|\f(\cdot,\o)\|_{L^\be(E)}\|\psi(\cdot,\o)\|_{L^\be(F)}\|\chi(\cdot,\o)\|_{L^\be(G)}\,d\s(\o)<\be.
$$

Suppose now that $T_1$ is a bounded linear operator from $\h_2$ to $\h_1$ and $T_2$ is a bounded linear operator from $\h_3$ to $\h_2$. For a function $\Phi$ in
$L^\be(E_1)\hat\otimes_{\rm i}L^\be(E_2)\hat\otimes_{\rm i}L^\be(E_3)$ of the form \rf{ttp}, we put
\begin{align}
\label{opr}
&\int\limits_\X\int\limits_\Y\int\limits_\cZ\Phi(x,y,z)
\,d E_1(x)T_1\,dE_2(y)T_2\,dE_3(z)\\[.2cm]
\df&\int\limits_\O\left(\,\int\limits_\X \f(x,\o)\,dE_1(x)\right)T_1
\left(\,\int\limits_\Y \psi(y,\o)\,dE_2(y)\right)T_2
\left(\,\int\limits_\cZ \chi(z,\o)\,dE_3(z)\right)\,d\s(\o).\nonumber
\end{align}

It was shown in \cite{Pe5} (see also \cite{ACDS} for a different proof)  that the above definition does not depend on the choice of a representation \rf{ttp}.
%

It is easy to see that the following inequality holds
$$
\left\|\int\limits_\X\int\limits_\Y\int\limits_\cZ\Phi(x,y,z)
\,dE_1(x)T_1\,dE_2(y)T_2\,dE_3(z)\right\|
\le\|\Phi\|_{L^\be\hat\otimes_{\rm i}L^\be\hat\otimes_{\rm i}L^\be}\cdot\|T_1\|\cdot\|T_2\|.
$$

In particular, the triple operator integral on the left-hand side of \rf{opr} can be defined if $\Phi$ belongs to the projective
tensor product $L^\be(E_1)\hat\otimes L^\be(E_2)\hat\otimes L^\be(E_3)$, i.e., $\Phi$ admits a representation
$$
\Phi(x,y,z)=\sum_{n\ge1}\f_n(x)\psi_n(y)\chi_n(z),
$$
where $\f_n\in L^\be(E_1)$, $\psi_n\in L^\be(E_2)$, $\chi_n\in L^\be(E_3)$ and
$$
\sum_{n\ge1}\|\f_n\|_{L^\be(E_1)}\|\psi_n\|_{L^\be(E_2)}\|\chi_n\|_{L^\be(E_3)}<\be.
$$

In a similar way one can define multiple operator integrals, see \cite{Pe5}.

Recall that multiple operator integrals were considered earlier in \cite{Pa} and \cite{S}. However, in those papers the class of functions 
$\Phi$ for which the left-hand side of \rf{opr} was defined is much narrower than in the definition given above.

Multiple operator integrals are used in \cite{Pe5} in connection with the problem of evaluating higher order operator derivatives. 
To obtain formulae for higher operator derivatives, one has to integrate divided differences of higher orders (see \cite{Pe5}). 

In this paper we are going to integrate divided differences of higher orders to estimate the norms of higher order differences \rf{hod}.

For a function $f$ on the circle the  {\it divided differences $\dg^k f$ of order $k$} are defined inductively as follows:
$$
\dg^0f\df f;
$$
if $k\ge1$, then in the case when $\l_1,\l_2,\cdots,\l_{k+1}$ are distinct points in $\T$,
$$
(\dg^{k}f)(\l_1,\cdots,\l_{k+1})\df
\frac{(\dg^{k-1}f)(\l_1,\cdots,\l_{k-1},\l_k)-
(\dg^{k-1}f)(\l_1,\cdots,\l_{k-1},\l_{k+1})}{\l_{k}-\l_{k+1}}
$$
(the definition does not depend on the order of the variables). Clearly,
$$
\dg f=\dg^1f.
$$
If $f\in C^k(\T)$, then $\dg^{k}f$ extends by continuity to a function defined for all points $\l_1,\l_2,\cdots,\l_{k+1}$.

It can be shown that
$$
({\frak D}^n f)(\l_1,\dots,\l_{n+1})=\sum\limits_{k=1}^{n+1}f(\l_k)
\prod\limits_{j=1}^{k-1}(\l_k-\l_j)^{-1}\prod\limits_{j=k+1}^{n+1}(\l_k-\l_j)^{-1}.
$$

Similarly, one can define the divided difference of order $k$ for functions on the real line.

It was shown in \cite{Pe5} that if $f$ is a trigonometric polynomial of degree $d$, then
\bay
\label{Bok}
\big\|\dg^k f\big\|_{C(\T)\hat\otimes\cdots\hat\otimes C(\T)}\le\const d^k\|f\|_{L^\be}.
\ey
If $f$ is an entire function of exponential type at most $\s$ that is bounded on $\R$, then
\bay
\label{Boke}
\big\|\dg^k f\big\|_{L^\be\hat\otimes_{\rm i}\cdots\hat\otimes_{\rm i} L^\be}\le\const \s^k\|f\|_{L^\be(\R)}.
\ey

Note that recently in \cite{JTT} Haagerup tensor products were used to define multiple operator integrals. However, it is
not clear whether this can lead to stronger results in perturbation theory.

\medskip

{\bf 3.3. Multiple operator integrals with respect to semi-spectral measures.}
Let $\h$ be a Hilbert space and let $(\X,{\frak B})$ be a measurable space.
A map $\E$ from ${\frak B}$ to the algebra ${\mathscr B}(\h)$ of all bounded operators on $\h$ is called a {\it semi-spectral measure}
if 
$$
\E(\D)\ge\0,\quad\D\in{\frak B},
$$
$$
\E(\varnothing)=\0\quad\mbox{and}\quad\E(\X)=I,
$$
and for a sequence $\{\D_j\}_{j\ge1}$ of disjoint sets in ${\frak B}$,
$$
\E\left(\bigcup_{j=1}^\be\D_j\right)=\lim_{N\to\be}\sum_{j=1}^N\E(\D_j)\quad\mbox{in the weak operator topology}.
$$

\medskip

If $\K$ is a Hilbert space, $(\X,{\frak B})$ is a measurable space,  $E:{\frak B}\to{\mathscr B}(\K)$ is a spectral measure, and $\h$ is
a subspace of $\K$, then it is easy to see that the map $\E:{\frak B}\to {\mathscr B}(\h)$ defined by
\bay
\label{dil}
\E(\D)=P_\h E(\D)\big|\h,\quad\D\in{\frak B},
\ey
is a semi-spectral measure. Here $P_\h$ stands for the orthogonal projection onto $\h$.

Naimark proved in \cite{N}  that all semi-spectral measures can be obtained in this way, i.e.,
a semi-spectral measure is always a {\it compression} of a spectral measure. A spectral measure $E$ satisfying \rf{dil} is called a {\it spectral dilation of the semi-spectral measure} $\E$.

A spectral dilation $E$ of a semi-spectral measure $\E$ is called {\it minimal} if 
$$
\K=\clos\spn\{E(\D)\h:~\D\in{\frak B}\}.
$$

It was shown in \cite{MM} that if $E$ is a minimal spectral dilation of a semi-spectral measure $\E$, then
$E$ and $\E$ are mutually absolutely continuous and all minimal spectral dilations of a semi-spectral measure are isomorphic in the natural sense.

If $\f$ is a bounded complex-valued measurable function on $\X$ and $\E:{\frak B}\to {\mathscr B}(\h)$ is a semi-spectral measure, then the integral
\bay
\label{iss}
\int_\X \f(x)\,d\E(x)
\ey
can be defined as
\bay
\label{voi}
\int_\X \f(x)\,d\E(x)=\left.P_\h\left(\int_\X \f(x)\,d E(x)\right)\right|\h,
\ey
where $E$ is a spectral dilation of $\E$. It is easy to see that the right-hand side of \rf{voi} does not depend on the choice
of a spectral dilation. The integral \rf{iss} can also be computed as the limit of sums
$$
\sum \f(x_\a)\E(\D_\a),\quad x_\a\in\D_\a,
$$
over all finite measurable partitions $\{\D_\a\}_\a$ of $\X$.

If $T$ is a contraction on a Hilbert space $\h$, then by the Sz.-Nagy dilation theorem
(see \cite{SNF}),  $T$ has a unitary dilation, i.e., there exist a Hilbert space $\K$ such that
$\h\subset\K$ and a unitary operator $U$ on $\K$ such that
\bay
\label{DT}
T^n=P_\h U^n\big|\h,\quad n\ge0,
\ey
where $P_\h$ is the orthogonal projection onto $\h$. Let $E_U$ be the spectral measure of $U$.
Consider the operator set function $\E$ defined on the Borel subsets of the unit circle $\T$ by
$$
\E(\D)=P_\h E_U(\D)\big|\h,\quad\D\subset\T.
$$
Then $\E$ is a semi-spectral measure. It follows immediately from
\rf{DT} that 
\bay
\label{step}
T^n=\int_\T \z^n\,d\E(\z)=P_\h\int_\T\z^n\,dE_U(\z)\Big|\h,\quad n\ge0.
\ey
Such a semi-spectral measure $\E$ is called a {\it semi-spectral measure} of $\T$.
Note that it is not unique. To have uniqueness, we can consider a minimal unitary dilation $U$ of $T$,
which is unique up to an isomorphism (see \cite{SNF}).

It follows easily from  \rf{step} that 
$$
f(T)=P_\h\int_\T f(\z)\,dE_U(\z)\Big|\h
$$
for an arbitrary function $\f$ in the disk-algebra $C_A$.

In \cite{Pe2} and  \cite{Pe6} double operator integrals and multiple operator integrals with respect to semi-spectral measures were introduced.

Suppose that $(\X_1,{\frak B}_1)$ and $(\X_2,{\frak B}_2)$ are measurable spaces, and
$\E_1:{\frak B}_1\to{\mathscr B}(\h_1)$ and $\E_2:{\frak B}_2\to {\mathscr B}(\h_2)$ are semi-spectral measures.
Then double operator integrals
$$
\iint\limits_{\X_1\times\X_2}\Phi(x_1,x_2)\,d\E_1(x_1)Q\,d\E_2(X_2).
$$
were defined in \cite{Pe6} in the case when $Q\in\bS_2$ and $\Phi$ is a bounded measurable function. Double operator integrals were also defined in \cite{Pe6} in the case when $Q$ is a bounded linear operator and $\Phi$ belongs to the integral projective tensor product of the spaces $L^\be(\E_1)$
and $L^\be(\E_2)$.

In particular, the following analog of the Birman--Solomyak formula holds:
\bay
\label{BSc}
f(R)-f(T)=\iint\limits_{\T\times\T}\big(\dg f\big)(\z,\t)\,d\E_R(\z)(R-T)\,d\E_T(\t).
\ey
Here $T$ and $R$ contractions on Hilbert space, $\E_T$ and $\E_R$ are their semi-spectral measures, and $f$ is an analytic  function in $\dd$ of class
$\big(B_{\be1}^1\big)_+$.

Similarly, multiple operator integrals with respect to semi-spectral measures were defined in \cite{Pe6} for functions that belong to the integral projective tensor product of the corresponding $L^\be$ spaces.

We also mention here the paper \cite{KS}, in which another approach is used to study perturbations of functions of contractions.

\

\section{\bf H\"older--Zygmund estimates for self-adjoint operators}
\setcounter{equation}{0}
\label{s4}

\

In this section we show that H\"older functions on $\R$ of order $\a$, $0<\a<1$, must also be operator H\"older of order $\a$. We also obtain similar results for all  H\"older--Zygmund classes $\L_\a(\R)$, $0<\a<\be$. For simplicity, we give complete proofs in the case of bounded self-adjoint operators and 
explain without details that similar inequalities also hold for unbounded self-adjoint operators. We are going to give a detailed treatment of the case of unbounded operators in \cite{AP3}. 

We compare in this section our results with an inequality by Birman, Koplienko, and Solomyak \cite{BKS}.

Note that if $A$ and $B$ are self-adjoint operators, we say that the operator $A-B$ is bounded if $B=A+K$ for some bounded
self-adjoint operator $K$. In particular, this implies that the domains of $A$ and $B$ coincide. We say that $\|A-B\|=\be$ if there is no such a bounded operator $K$ that $B=A+K$.


\begin{thm}
\label{saH}
Let $0<\a<1$. Then there is a constant $c>0$ such that for 
every $f\in\L_\a(\R)$ and for arbitrary self-adjoint operators $A$ and $B$ on Hilbert space the following inequality holds:
$$
\|f(A)-f(B)\|\le c\,\|f\|_{\L_\a(\R)}\cdot\|A-B\|^\a.
$$
\end{thm}

\Pf If $A$ and $B$ are bounded operators, it follows from Theorem \ref{prodol} that we may assume that $f\in L^\be(\R)$ and we have to obtain an estimate for $\|f(A)-f(B)\|$ that does not depend on $\|f\|_{L^\be}$. 

Put
$$
f_n=f*W_n+f*W_n^\sharp.
$$
Let us show that
\bay
\label{skho}
f(A)-f(B)=\sum_{n=-\be}^\be\big(f_n(A)-f_n(B)\big)
\ey
and the series on the right converges absolutely in the operator norm.

For $N\in\Z$, we put
$$
g_N=f*V_N
$$
(recall that $V_N$ is the de la Vall\'ee Possin type kernel defined in \S\,2.1). Clearly,
$$
f=f*V_N+\sum_{n>N} f_n
$$
and the series on the right converges absolutely in the $L^\be$ norm.
Thus
$$
f(A)=\big(f*V_N\big)(A)+\sum_{n>N} f_n(A)\quad\mbox{and}\quad
f(B)=\big(f*V_N\big)(B)+\sum_{n>N} f_n(B),
$$
and the series converge absolutely in the operator norm. We have
\begin{align*}
f(A)-f(B)-\sum_{n>N}\big(f_n(A)-f_n(B)\big)&=\left(f(A)-
\sum_{n>N} f_n(A)\right)-\left(f(B)-\sum_{n>N} f_n(B)\right)\\[.2cm]
&=g_N(A)-g_N(B).
\end{align*}
Since $g_N\in L^\be(\R)$ and $g_N$ is an entire function of exponential type at most $2^{N+1}$,
it follows from \rf{BSF}  and \rf{Be} that
$$
\|g_N(A)-g_N(B)\|\le\const2^N\|f*V_N\|_{L^\be}\|A-B\|\le\const2^N\|f\|_{L^\be}\|A-B\|\to0
$$
as $N\to-\be$. This proves \rf{skho}.

Let now  $N$ be the integer such that
\bay
\label{AB}
2^{-N}<\|A-B\|\le2^{-N+1}.
\ey

We have 
$$
f(A)-f(B)=\sum_{n\le N}\big(f_n(A)-f_n(B)\big)+\sum_{n>N}\big(f_n(A)-f_n(B)\big).
$$

It follows from \rf{Wn} and \rf{AB} that
\begin{align*}
\left\|\sum_{n\le N}\big(f_n(A)-f_n(B)\big)\right\|&\le\sum_{n\le N}\big\|f_n(A)-f_n(B)\big\|\\[.2cm]
&\le\const\sum_{n\le N}2^n\|f_n\|_{L^\be}\|A-B\|\\[.2cm]
&\le\sum_{n\le N}2^n2^{-n\a}\|f\|_{\L_\a(\R)}\|A-B\|\\[.2cm]
&\le\const2^{N(1-\a)}\|f\|_{\L_\a(\R)}\|A-B\|\le\|f\|_{\L_\a(\R)}\|A-B\|^\a.
\end{align*}

On the other hand, 
\begin{align*}
\left\|\sum_{n>N}\big(f_n(A)-f_n(B)\big)\right\|
&\le\sum_{n>N}\Big(\|f_n(A)\|+\|f_n(B)\|\Big)\\[.2cm]
&\le2\sum_{n>N}\|f_n\|_{L^\be}\le\const\sum_{n>N}2^{-N\a}\|f\|_{\L_\a(\R)}\\[.2cm]
&\le\const2^{-N\a}\|f\|_{\L_\a(\R)}\le\const\|f\|_{\L_\a(\R)}\|A-B\|^\a
\end{align*}
by \rf{AB}. This completes the proof in the case of bounded self-adjoint operators. 

In the case of unbounded self-adjoint operators the same reasoning holds if by \lb$f(A)-f(B)$ we understand the series
$$
\sum_{n\in\Z}\big(f_n(A)-f_n(B)\big),
$$
which converges absolutely. In \cite{AP3} we are going to consider the case of unbounded self-adjoint operators in more detailed.
$\bl$

\medskip

{\bf Remark.} Note that Birman, Koplienko, and Solomyak obtained in \cite{BKS} the following result: if $A$ and $B$ are positive self-adjoint operators and $0<\a<1$, then
$$
\|A^\a-B^\a\|\le\|A-B\|^\a.
$$
It follows from Theorem \ref{saH} that under the same assumptions 
$$
\|A^\a-B^\a\|\le\const\|A-B\|^\a.
$$
Indeed, it suffices to apply Theorem \ref{saH} to the operators $A$, $B$ and the function $f\in\L_\a(\R)$ defined by
$f(t)=|t|^\a$, $t\in\R$.

\medskip

Let us now state the result for arbitrary H\"older--Zygmund classes $\L_\a(\R)$. 

\begin{thm}
\label{sam}
Let $0<\a<m$. Then there exists a constant $c>0$ such that for every self-adjoint operators  $A$ and $K$ on Hilbert space the following inequality holds:
$$
\left\|\sum_{j=0}^m(-1)^{m-j}\left(\begin{matrix}m\\j\end{matrix}\right)f\big(A+jK\big)\right\|
\le c\,\|f\|_{\L_\a(\R)}\cdot\|K\|^\a.
$$
\end{thm}

We need the following lemma.

\begin{lem}
\label{m}
Let $m$ be a positive integer and let $f$ be a bounded function of class $B^m_{\be1}(\R)$. 
If $A$ and $K$ are self-adjoint operators
on Hilbert space, then
\begin{align*}
\sum_{j=0}^m&(-1)^{m-j}\left(\begin{matrix}m\\j\end{matrix}\right)f\big(A+jK\big)\\[.2cm]
&=
m!\underbrace{\int\cdots\int}_{m+1}(\dg^{m}f)(x_1,\cdots,x_{m+1})
\,dE_A(x_1)K\,dE_{A+K}(x_2)K\cdots K\,dE_{A+mK}(x_{m+1}).
\end{align*}
\end{lem}

For simplicity, we prove Theorem \ref{sam} and Lemma \ref{m} for $m=2$. The general case can be treated in the same way.

\medskip

{\bf Proof of Lemma \ref{m}.} In the case $m=2$ we have to establish the following formula for $f\in B_{\be1}^2(\R)$:
$$
f(A+K)-2f(A)+f(A-K)=\!2\!\iiint\!(\dg^2f)(x,y,z)\,dE_{A+K}(x)K\,dE_A(y)K\,dE_{A-K}(z).
$$
Put $T=f(A+K)-2f(A)+f(A-K)$. By \rf{BSF},

\begin{align*}
T&=
f(A+K)-f(A)-\big(f(A)-f(A-K)\big)\\[.2cm]
&=\iint(\dg f)(x,y)\,dE_{A+K}(x)K\,dE_A(y)-
\iint(\dg f)(x,y)\,dE_{A}(x)K\,dE_{A-K}(y)\\[.2cm]
&=\iint(\dg f)(x,y)\,dE_{A+K}(x)K\,dE_A(y)-
\iint(\dg f)(x,y)\,dE_{A+K}(x)K\,dE_{A-K}(y)\\[.2cm]
&+\iint(\dg f)(x,y)\,dE_{A+K}(x)K\,dE_{A-K}(y)-
\iint(\dg f)(x,y)\,dE_{A}(x)K\,dE_{A-K}(y).
\end{align*}

We have
\begin{align*}
&\iint(\dg f)(x,y)\,dE_{A+K}(x)K\,dE_A(y)-
\iint(\dg f)(x,y)\,dE_{A+K}(x)K\,dE_{A-K}(y)\\[.2cm]
=&\iint(\dg f)(x,y)\,dE_{A+K}(x)K\,dE_A(y)-
\iint(\dg f)(x,z)\,dE_{A+K}(x)K\,dE_{A-K}(z)\\[.2cm]
=&\iiint(\dg f)(x,y)\,dE_{A+K}(x)K\,dE_A(y)\,dE_{A-K}(z)\\[.2cm]
-&\iiint(\dg f)(x,z)\,dE_{A+K}(x)K\,dE_A(y)\,dE_{A-K}(z)\\[.2cm]
=&\iiint(y-z)(\dg^2f)(x,y,z)\,dE_{A+K}(x)K\,dE_A(y)\,dE_{A-K}(z)\\[.2cm]
=&\iiint(\dg^2f)(x,y,z)\,dE_{A+K}(x)K\,dE_A(y)K\,dE_{A-K}(z).
\end{align*}

Similarly,
\begin{align*}
&\iint(\dg f)(x,y)\,dE_{A+K}(x)K\,dE_{A-K}(y)-
\iint(\dg f)(x,y)\,dE_{A}(x)K\,dE_{A-K}(y)\\[.2cm]
=&\iiint(\dg^2f)(x,y,z)\,dE_{A+K}(x)K\,dE_A(y)K\,dE_{A-K}(z).
\end{align*}

Thus
$$
T=2\iiint(\dg^2f)(x,y,z)\,dE_{A+K}(x)K\,dE_A(y)K\,dE_{A-K}(z).\quad\bl
$$

\medskip

{\bf Proof of Theorem \ref{sam}.}
By Theorem \ref{prodol}, we may assume that $f$ is a bounded function.

We are going to use the same notation $f_n$ and $g_N$ as in the proof of Theorem \ref{saH}. In the case when $A$ and $K$ are bounded
self-adjoint operators we show that 
\bay
\label{skh}
f(A+K)-2f(A)+f(A-K)=\sum_{n=-\be}^\be\big(f_n(A+K)-2f_n(A)+f_n(A-K)\big),
\ey
and the series converges absolutely in the operator norm. As in the proof of Theorem \ref{saH}, we can easily see that
$$
f(A+K)=\big(f*V_N\big)(A+K)+\sum_{n>N} f_n(A+K),
$$
$$
f(A)=\big(f*V_N\big)(A)+\sum_{n>N} f_n(A),
$$
and
$$
f(A-K)=\big(f*V_N\big)(A-K)+\sum_{n>N} f_n(A-K),
$$
and the series converge absolutely in the operator norm.
It follows that
\begin{align*}
f(A+K)&-2f(A)+f(A-K)-\sum_{n>N}\big(f_n(A+K)-2f_n(A)+f_n(A-K)\big)\\[.2cm]
&=\left(f(A+K)-
\sum_{n>N} f_n(A+K)\right)-2\left(f(A)-\sum_{n>N} f_n(A)\right)\\[.2cm]
&+\left(f(A-K)-\sum_{n>N} f_n(A-K)\right)
=g_N(A+K)-2g_N(A)+g_N(A-K).
\end{align*}
Since $g_N\in L^\be(\R)$ and $g_N$ is an entire function of exponential type at most $2^{N+1}$,
it follows from Lemma \ref{m} and from \rf{Boke} that
\begin{align*}
\big\|g_N(A+K)-2g_N(A)+g_N(A-K)\big\|&\le\const2^{-2N}\|g_N\|_{L^\be}\|K\|\\[.2cm]
&\le\const2^{-2N}\|f\|_{L^\be}\|K\|\to0\quad\mbox{and}\quad m\to\be.
\end{align*}
This implies that the series on the right-hand side of \rf{skh} converges absolutely in the operator norm.

As in the proof of Theorem \ref{saH}, we consider the integer $N$ satisfying
\bay
\label{K}
2^{-N}<\|K\|\le2^{-N+1}.
\ey
Put now
$$
T_1\df\sum_{n\le N}\big(f_n(A+K)-2f_n(A)+f_n(A-K)\big)
$$
and
$$
T_2\df\sum_{n>N}\big(f_n(A+K)-2f_n(A)+f_n(A-K)\big).
$$

It follows now from Lemma \ref{m}, from \rf{K},  and \rf{Boke} that
\begin{align*}
\|T_1\|&\le
\sum_{n\le N}\|f_n(A+K)-2f_n(A)+f_n(A-K)\|\\[.2cm]
&=2\sum_{n\le N}\left\|\iiint(\dg^2f_n)(x,y,z)\,dE_{A+K}(x)K\,dE_A(y)K\,dE_{A-K}(z)\right\|\\[.2cm]
&\le\const\sum_{n\le N}2^{2n}\|f_n\|_{L^\be}\|K\|^2\le\const\sum_{n\le N}2^{n(2-\a)}\|f\|_{\L_\a(\R)}\|K\|^2\\[.2cm]
&\le\const2^{N(2-\a)}\|K\|^2\|f\|_{\L_\a(\R)}\le\const\|f\|_{\L_\a(\R)}\|K\|^\a.
\end{align*}
On the other hand, by \rf{K},
\begin{align*}
\|T_2\|&\le\sum_{n>N}\big\|\big(f_n(A+K)-2f_n(A)+f_n(A-K)\big)\big\|\\[.2cm]
&\le4\sum_{n>N}\|f_n\|_{L^\be}\le\const\sum_{n>N}2^{-n\a}\|f\|_{\L_\a(\R)}\\[.2cm]
&\le\const2^{-N\a}\|f\|_{\L_\a(\R)}\le\const\|K\|^\a.
\end{align*}

As in the case $\a<1$, for unbounded self-adjoint operators we understand by \lb$f(A+K)-2f(A)+f(A-K)$ the sum of the following series
$$
\sum_{n\in\Z}\big(f_n(A+K)-2f_n(A)+f_n(A-K)\big),
$$
which converges absolutely. We refer the reader to \cite{AP3} where the case of unbounded self-adjoint operators will be considered in more detail. $\bl$

\begin{cor}
\label{Z}
There exists a constant $c>0$ such that  for an arbitrary function $f$ in the Zygmund class $\L_1(\R)$ and arbitrary self-adjoint operators $A$ and $K$, the following inequality holds:
$$
\big\|f(A+K)-2f(A)+f(A-K)\big\|\le c\|f\|_{\L_1(\R)}\|K\|.
$$
\end{cor}

{\bf Remark.} We can interpret Theorem \ref{sam} in the following way. Consider the measure $\nu$ on $\R$ defined by
$$
\nu\df\D_1^m\d_0=\sum_{j=0}^m(-1)^{m-j}\left(\begin{matrix}m\\j\end{matrix}\right)\d_{-j},
$$
where for $a\in\R$, $\d_a$ is the unit point mass at $a$. Then
$$
\sum_{j=0}^m(-1)^{m-j}\left(\begin{matrix}m\\j\end{matrix}\right)f\big(A+jK\big)=
\int_\R f(A-tK)\,d\nu(t).
$$
Clearly, $\nu$ determines a continuous linear functional on $\l_\a(\R)$ defined by
$$
f\mapsto \int_\R f(t)\,d\nu(t).
$$
In other words, $\nu\in B_1^{-\a}(\R)$ (see \S\,2.1).
We are going to generalize Theorem \ref{sam} to the case of an arbitrary distribution in $B_1^{-\a}(\R)$. 

For simplicity, we consider here the case of bounded self-adjoint operators $A$. In \cite{AP3}
we will consider the case of an arbitrary (not necessarily bounded) self-adjoint operator $A$.

It follows from Theorem \ref{sam} that for arbitrary vectors $u$ and $v$ in our Hilbert space $\h$ and for an arbitrary function $f$ in $\L_\a(\R)$, the function
$$
t\mapsto f_{A,K}^{u,v}(t)\df \big(f(A-tK)u,v\big)
$$
belongs to $\L_\a(\R)$. Identifying the space $\L_\a(\R)$ with the dual space to $B_1^{-\a}$ (see \S\,2.1), we can consider
for every distribution $g$ in $B_1^{-\a}(\R)$ the value $\langle f_{A,K}^{u,v},g\rangle$ of $f_{A,K}^{u,v}\in \big(B_1^{-\a}(\R)\big)^*$ at $g$.
We define now the operator $\Q_{A,K}^g:\L_\a(\R)\to{\mathscr B}(\h)$ by
$$
\Big(\big(\Q_{A,K}^gf\big)u,v\Big)=\langle f_{A,K}^{u,v},g\rangle, \quad f\in \L_\a(\R),\quad u,\,v\in\h.
$$

\begin{thm}
\label{QAKg}
Let $\a>0$. Then there exists $c>0$ such that for every self-adjoint operators $A$ and $K$, for every $f\in\L_\a(\R)$, and for every $g\in B_1^{-\a}$,
\bay
\label{QAKgf}
\big\|\Q_{A,K}^gf\|\le c\,\|f\|_{\L_\a(\R)}\|g\|_{B_1^{-\a}(\R)}\|K\|^\a.
\ey
\end{thm}

\Pf Let $m$ be the smallest integer greater than $\a$.
By Theorem \ref{sam}, inequality \rf{QAKgf} holds for $g=\D_1^m\d_0$. Hence, the result also holds for
$g=\D_h^m\d_a$ for arbitrary $h,\,a\in\R$. 

To complete the proof, it suffices to use the following fact (see \cite{A}, Th. 3.1): if $g\in B_1^{-\a}(\R)$, then $g$ admits a representation
in the form of a norm convergent series
$$
g=\sum_{j\ge1}\l_j\D^m_{h_j}\d_{a_j},\quad h_j,\,a_j\in\R,
$$
such that 
$$
\sum_{j\ge1}|\l_j|\cdot\big\|\D^m_{h_j}\d_{a_j}\big\|_{B_1^{-\a}(\R)}\le\const\|g\|_{B_1^{-\a}(\R)}.\quad\bl
$$

\

\section{\bf The case of unitary operators}
\setcounter{equation}{0}
\label{s5}

\

In this section we obtain analogs of the results of the previous section for functions of unitary operators. We also obtain an estimate
for $\|f(U)-f(V)\|$ for a function $f$ in the Zygmund class $\L_1$ and unitary operators $U$ and $V$.

\begin{thm}
\label{uH}
Let $0<\a<1$. Then there is a constant $c>0$ such that for 
every $f\in\L_\a$ and for arbitrary unitary operators $U$ and $V$ on Hilbert space the following inequality holds:
$$
\|f(U)-f(V)\|\le c\,\|f\|_{\L_\a}\cdot\|U-V\|^\a.
$$
\end{thm}

\Pf Let $f\in\L_\a$. We have
$$
f=\pp_+f+\pp_-f=f_++f_-.
$$
We estimate $\|f_+(U)-f_+(V)\|$. The norm of  $f_-(U)-f_-(V)$ can be obtained in the same way.
Thus we assume that $f=f_+$. Let
$$
f_n\df f*W_n.
$$
Then
\bay
\label{dia}
f=\sum_{n\ge0}f_n.
\ey
Clearly, we may assume that $U\ne V$. Let $N$ be the nonnegative integer such that
\bay
\label{N}
2^{-N}<\|U-V\|\le2^{-N+1}.
\ey
We have 
$$
f(U)-f(V)=\sum_{n\le N}\big(f_n(U)-f_n(V)\big)+\sum_{n> N}\big(f_n(U)-f_n(V)\big).
$$
By the Birman--Solomyak formula for unitary operators and by \rf{Bp},
\begin{align*}
\left\|\sum_{n\le N}\big(f_n(U)-f_n(V)\big)\right\|&\le\sum_{n\le N}\big\|f_n(U)-f_n(V)\big\|\\[.2cm]
&\le\const\sum_{n\le N}2^n\|U-V\|\cdot\|f_n\|_{L^\be}\\[.2cm]
&\le\const\|U-V\|\sum_{n\le N}2^n2^{-n\a}\|f\|_{\L_\a}\\[.2cm]
&\le\const\|U-V\|2^{N(1-\a)}\|f\|_{\L_\a}\le\const\|U-V\|^\a\|f\|_{\L_\a},
\end{align*}
the last inequality being a consequence of \rf{N}.

On the other hand,
\begin{align*}
\left\|\sum_{n>N}\big(f_n(U)-f_n(V)\big)\right\|&\le\sum_{n>N}2\|f_n\|_{L^\be}
\le\const\sum_{n>N}2^{-n\a}\|f\|_{\L_\a}\\[.2cm]
&\le\const2^{-N\a}\|f\|_{\L_\a}\le\const\|U-V\|^\a\|f\|_{\L_\a}.\quad\bl
\end{align*}

To obtain an analog of Theorem \ref{sam} for unitary operator, we are going to represent a finite difference
$$
\sum_{j=1}^N(-1)^{j-1}\left(\begin{matrix}N-1\\j-1\end{matrix}\right)f(U_j)
$$
for unitary operators $U_1,\cdots,U_N$ as a linear combination of multiple operator integrals. 

Note that algebraic formulae in the case of unitary operators are more complicated than in the case of self-adjoint operators.
That is why we consider the case of unitary operators in more detail.

We first illustrate the idea in the special case $N=3$. 

Let us show that for unitary operators $U_1$, $U_2$ and $U_3$ and for $f\in B_{\be1}^2$, 
\begin{align}
\label{N3}
f(U_1)-2f(U_2)+f(U_3)&=\!2\!\iiint\!(\dg^2f)(\z,\t,\up)\,dE_1(\z)(U_1-U_2)\,dE_2(\t)(U_2-U_3)\,dE_3(\up)\nonumber\\[.2cm]
&+\iint(\dg f)(\z,\t)\,dE_1(\z)(U_1-2U_2+U_3)\,dE_3(\t).
\end{align}

Indeed, let $T=f(U_1)-2f(U_2)+f(U_3)$. Then
\begin{align*}
T&=
f(U_1)-f(U_2)-\big(f(U_2)-f(U_3)\big)\\[.2cm]
&=\iint(\dg f)(\z,\t)\,dE_1(\z)(U_1-U_2)\,dE_2(\t)-
\iint(\dg f)(\z,\t)\,dE_2(\z)(U_2-U_3)\,dE_3(\t)\\[.2cm]
&=\iint(\dg f)(\z,\t)\,dE_1(\z)(U_1-U_2)\,dE_2(\t)-
\iint(\dg f)(\z,\t)\,dE_1(\z)(U_1-U_2)\,dE_3(\t)\\[.2cm]
&+\iint(\dg f)(\z,\t)\,dE_1(\z)(U_1-U_2)\,dE_3(\t)
-\iint(\dg f)(\z,\t)\,dE_1(\z)(U_2-U_3)\,dE_3(\t)
\\[.2cm]
&+\iint(\dg f)(\z,\t)\,dE_1(\z)(U_2-U_3)\,dE_3(\t)
-\iint(\dg f)(\z,\t)\,dE_2(\z)(U_2-U_3)\,dE_3(\t).
\end{align*}

We have
\begin{align*}
&\iint(\dg f)(\z,\t)\,dE_1(\z)(U_1-U_2)\,dE_2(\t)-
\iint(\dg f)(\z,\t)\,dE_1(\z)(U_1-U_2)\,dE_3(\t)\\[.2cm]
=&\iint(\dg f)(\z,\t)\,dE_1(\z)(U_1-U_2)\,dE_2(\t)-
\iint(\dg f)(\z,\up)\,dE_1(\z)(U_1-U_2)\,dE_3(\up)\\[.2cm]
=&\iiint(\dg f)(\z,\t)\,dE_1(\z)(U_1-U_2)\,dE_2(\t)\,dE_3(\up)\\[.2cm]
-&\iiint(\dg f)(\z,\up)\,dE_1(\z)(U_1-U_2)\,dE_2(\t)\,dE_3(\up)\\[.2cm]
=&\iiint(\t-\up)(\dg^2f)(\z,\t,\up)\,dE_1(\z)(U_1-U_2)\,dE_2(\t)\,dE_3(\up)\\[.2cm]
=&\iiint(\dg^2f)(\z,\t,\up)\,dE_1(\z)(U_1-U_2)\,dE_2(\t)(U_2-U_3)\,dE_3(\up).
\end{align*}
Similarly,
\begin{align*}
&\iint(\dg f)(\z,\t)\,dE_1(\z)(U_2-U_3)\,dE_3(\t)
-\iint(\dg f)(\z,\t)\,dE_2(\z)(U_2-U_3)\,dE_3(\t)\\[.2cm]
=&\iiint(\dg^2f)(\z,\t,\up)\,dE_1(\z)(U_1-U_2)\,dE_2(\t)(U_2-U_3)\,dE_3(\up).
\end{align*}
Finally,
\begin{align*}
&\iint(\dg f)(\z,\t)\,dE_1(\z)(U_1-U_2)\,dE_3(\t)
-\iint(\dg f)(\z,\t)\,dE_1(\z)(U_2-U_3)\,dE_3(\t)\\[.2cm]
=&\iint(\dg f)(\z,\t)\,dE_1(\z)(U_1-2U_2+U_3)\,dE_3(\t).\quad\bl
\end{align*}

Consider now the general case. Suppose that $\cU=\{U_j\}_1^N$ is a finite family of unitary operators. Denote by $E_j$ the spectral measure of $E_j$. For $1\le j<k\le N$, we put
$$
T(j,k)=\sum_{s=0}^{k-j}(-1)^s\left(\begin{matrix}k-j\\s\end{matrix}\right)U_{j+s}.
$$
Note that
\bay
\label{Tjk}
T(j,k)-T(j+1,k+1)=T(j,k+1),\quad1\le j<k\le N-1.
\ey

Let $J$ be a nonempty subset of $\{1,2,\cdots,N\}$. We denote by $d=d_J$ the number of elements of $J$. Suppose that $J=\{j_1,j_2,\cdots,j_d\}$,
where $j_1<j_2<\cdots<j_d$. For $f\in B_{\be1}^{d-1}$, we put
$$
\I_J(\cU,f)\df\underbrace{\int\cdots\int}_d\big(\dg^{d-1}f\big)(\z_1,\cdots,\z_d)\,dE_{j_1}(\z_1)\prod_{s=2}^dT(j_{s-1},j_s)\,dE_{j_s}(\z_s).
$$
Though, we need the case, $d_J\ge2$, but we still can assume that $d_J=1$, in which case we put
$$
\I_J(\cU,f)\df\int f(\z)\,dE_j(\z),\quad\mbox{where}\quad J=\{j\}.
$$ 

We denote by $\frak A$ the collection of all finite subsets of the set of positive integers and by ${\frak A}_N$ the collection of all subsets $J\in{\frak A}$ such that the maximal element of $J$ is $N$.

If $J_1,\,J_2\in{\frak A}$, we say that $J_1$ is an {\it ancestor of} $J_2$ if $J_2$ can be partitioned in nonempty subsets $J_2^{\prime}$ and 
$J_2^{\prime\prime}$ such that $\max J_2^{\prime}<\min J_2^{\prime\prime}$ and $J_1=J_2^\prime\bigcup\big(J_2^{\prime\prime}-1\big)$ (by $\L-1$ we mean the left translate of a subset $\L$ of $\Z$ by 1).
Each such partition is called an {\it evidence of the fact that $J_1$ is an ancestor of $J_2$}. We denote by $\#(J_1,J_2)$ the number of such evidences and we put $\#(J_1,J_2)=0$ if $J_1$ is not an ancestor of $J_2$. Note that the property of being an ancestor is not transitive.

If $\#(J_1,J_2)\ge1$, then $\max J_2=1+\max J_1$ and $0\le d_{J_2}-d_{J_1}\le1$. It is also easy to see that if $d_{J_1}=d_{J_2}$, then
$\#(J_1,J_2)=1$.

Let us construct now the family $\vk_J$ of integers by induction. Put $\vk_{\{1\}}=1$. Suppose that the numbers $\vk_J$ are defined for 
$J\in{\frak A}_{N-1}$. Let $J\in{\frak A}_N$. Put
$$
\vk_J=\sum_{I\in{\frak A}_{N-1}}\#(I,J)\vk_I.
$$
Clearly, $\vk_J$ is a positive integer for every $J\in{\frak A}$. We leave for the reader the verification of the fact that for $\{j_1,j_2,\cdots,j_d\}\in{\frak A}$,
$$
\vk_J=\frac{(j_d-j_1)!}{\prod_{s=2}^d\big(j_s-j_{s-1}\big)!}.
$$

\begin{thm}
\label{gen}
Let $N$ be a positive integer and let $\cU=\{U_j\}_{j=1}^N$ be unitary operators on Hilbert space.
Suppose that $f$ is a function in the Besov space $B_{\be1}^{N-1}$. Then
$$
\sum_{j=1}^N(-1)^{j-1}\left(\begin{matrix}N-1\\j-1\end{matrix}\right)f(U_j)
=\sum_{J\in{\frak A}_N}\vk_J\I_J(\cU,f).
$$
\end{thm}

We need one more lemma. To state it, we introduce some more notation. For $J\in{\frak A}$, we denote by ${\frak L}(J)$ the collection of nonempty proper
subsets of $J$ such that 
$$
\max\L<\min(J\setminus\L).
$$
For $\L\in{\frak L}(J)$, we put
$$
\L^\circ_J\df J\setminus\L\quad\mbox{and}\quad\L^\bullet_J\df \L^\circ_J\cup\{\max\L\}.
$$
If $J$ is specified, we write $\L^\circ$ and $\L^\bullet$ instead of $\L^\circ_J$ and $\L^\bullet_J$.

\begin{lem}
\label{OJ}
Let $J\in{\frak A}_{N-1}$. Then
$$
\I_J(\cU,f)\!-\I_{J+1}(\cU,f)=\!\!\sum_{\L\in{\frak L}(J)}\I_{\L\cup(\L^\circ+1)}(\cU,f)
+\!\!\sum_{\L\in{\frak L}(J)}\I_{\L\cup(\L^\bullet+1)}(\cU,f)+\I_{J\cup\{N\}}(\cU,f).
$$
\end{lem}

\Pf The above identity can be verified straightforwardly if we observe that 
the multiple operator integral
$$
\underbrace{\int\cdots\int}_d\big(\dg^{d-1}f\big)(\z_1,\cdots,\z_d)\,dF_{1}(\z_1)\prod_{s=2}^d Q_{s-1}\,dF_{s}(\z_s)
$$
is a multilinear function in the operators $Q_s$ and use the following easily verifiable identity:
$$
\iint\big({\frak D}f\big)(\z_1,\z_2)\,dE_1(\z_1)(U_1-U_2)\,dE_2(\z_2)=\int f(\z)\,dE_1(\z)-\int f(\z)\,dE_2(\z).\quad\bl
$$

\medskip

{\bf Proof of Theorem \ref{gen}.} We argue by induction on $N$. For $N=1$, we have
$$
f(U_1)=\int f(\z_1)\,dE(\z_1).
$$
Suppose that the result holds for $N-1$ unitary operators. Put $\cU^-\df\{U_{j+1}\}_{j=1}^{N-1}$. We have
$$
\sum_{j=1}^{N-1}(-1)^{j-1}\left(\begin{matrix}N-2\\j-1\end{matrix}\right)f(U_j)
=\sum_{J\in{\frak A}_{N-1}}\vk_J\I_J(\cU,f)
$$
and
$$
\sum_{j=1}^{N-1}(-1)^{j-1}\left(\begin{matrix}N-2\\j-1\end{matrix}\right)f(U_{j+1})=\sum_{J\in{\frak A}_{N-1}}\vk_J\I_J(\cU^-,f)=
\sum_{J\in{\frak A}_{N-1}}\vk_J\I_{J+1}(\cU,f).
$$
It follows now from \rf{Tjk} and Lemma \ref{OJ} that
\begin{align*}
\sum_{j=1}^N&(-1)^{j-1}\left(\begin{matrix}N-1\\j-1\end{matrix}\right)f(U_j)=\sum_{J\in{\frak A}_{N-1}}\vk_J\big(\I_J(\cU,f)-\I_{J+1}(\cU,f)\big)\\[.2cm]
&=\sum_{J\in{\frak A}_{N-1}}\vk_J
\left(\sum_{\L\in{\frak L}(J)}\I_{\L\cup(\L^\circ+1)}(\cU,f)
+\sum_{\L\in{\frak L}(J)}\I_{\L\cup(\L^\bullet+1)}(\cU,f)+\I_{J\cup\{N\}}(\cU,f)
\right).
\end{align*}
It remains to observe that a set $J$ in ${\frak A}_{N-1}$ is an ancestor of a set $J_0$ in ${\frak A}_N$ if and only if 
$J_0=\L\cup(\L^\circ+1)$ for some $\L\in{\frak L}(J)$ or $J_0=\L\cup(\L^\bullet+1)$ for some $\L\in{\frak L}(J)$ or $J_0=J\cup\{N\}$. $\bl$

\begin{thm}
\label{hou}
Let $m$ be a positive integer and $0<\a<m$. Then there exists a constant $c>0$ such that for every $f\in\L_\a$ and
for an arbitrary unitary operator $U$ and an arbitrary bounded self-adjoint operator $A$ on Hilbert space the following inequality holds:
$$
\left\|\sum_{k=0}^m(-1)^k\left(\begin{matrix}m\\k\end{matrix}\right)f\big(e^{{\rm i}kA}U\big)
\right\|\le c\,\|f\|_{\L_\a}\|A\|^\a.
$$
\end{thm}

\Pf For simplicity we give a proof for $m=2$. The general case can be treated in the same way. We have to show that
for $0<\a<2$, there is a constant $c>0$ such that for 
every $f\in\L_\a$ and for arbitrary unitary operators $U$ and $\V$ on Hilbert space the following inequality holds:
$$
\|f(\V U)-2f(U)+f(\V^*U)\|\le c\|f\|_{\L_\a}\|I-\V\|^\a. 
$$

%

As in the proof of Theorem \ref{uH}, we assume that $f=f_+$ and consider the expansion
$$
f=\sum_{n\ge0}f_n.
$$
Let $N$ be the nonnegative integer such that
\bay
\label{N1}
2^{-N}<\|I-\V\|\le2^{-N+1}.
\ey
We have
\begin{align*}
f(\V U)-2f(U)+f(\V^*U)&=
\sum_{n\le N}\big(f_n(\V U)-2f_n(U)+f_n(\V^*U)\big)\\[.2cm]
&+\sum_{n> N}\big(f_n(\V U)-2f_n(U)+f_n(\V^*U)\big).
\end{align*}

Let $T_n=f_n(\V U)-2f_n(U)+f_n(\V^*U)$.
It follows from \rf{N3} that
\begin{align*}
T_n&=2\iiint(\dg^2f_n)(\z,\t,\up)\,dE_{\V U}(\z)U(\V-I)\,dE_U(\t)U(I-\V^*)\,dE_{\V^*U}(\up)\\[.2cm]
&+\iint(\dg f_n)(\z,\t)\,dE_{\V U}(\z)U(\V-2I+\V^*)\,dE_{\V^*U}(\t).
\end{align*}
By \rf{Bok}, we have
$$
\left\|\iiint\!(\dg^2f_n)(\z,\t,\up)\,dE_{\V U}(\z)U(\V-I)\,dE_U(\t)U(I-\V^*)\,dE_{\V^*U}(\up)\right\|\!\le\!\const2^{2n}\|I-\V\|^2.
$$
On the other hand, by \rf{Bp},
\begin{align*}
\left\|\iint(\dg f_n)(\z,\t)\,dE_{\V U}(\z)U(\V-2I+\V^*)\,dE_{\V^*U}(\t)\right\|&\le\const2^n\|\V-2I+\V^*\|\\[.2cm]
&\le\const2^n\|I-\V\|^2.
\end{align*}
Thus
\begin{align*}
\left\|\sum_{n\le N}\big(f_n(\V U)-2f_n(U)+f_n(\V^*U)\big)\right\|&\le\const\|I-\V\|^2\sum_{n\le N}2^{2n}\|f_n\|_{L^\be}\\[.2cm]
&\le\const\|I-\V\|^2\sum_{n\le N}2^{2n}2^{-n\a}\|f\|_{\L_\a}\\[.2cm]
&\le\const\|I-\V\|^22^{N(2-\a)}\|f\|_{\L_\a}\\[.2cm]
&\le\const\|f\|_{\L_\a}\|I-\V\|^\a
\end{align*}
by \rf{N1}.

To complete the proof, we observe that
\begin{align*}
\left\|\sum_{n> N}\big(f_n(\V U)-2f_n(U)+f_n(\V^*U)\big)\right\|&\le
\sum_{n>N}\|(f_n(\V U)-2f_n(U)+f_n(\V^*U)\big)\|\\[.2cm]
&\le\sum_{n>N}4\|f_n\|_{L^\be}\le\const\|f\|_{\L_\a}\sum_{n>N}2^{-n\a}\\[.2cm]
&\le\const\|f\|_{\L_\a}2^{-N\a}\le\const\|I-\V\|^\a
\end{align*}
by \rf{N1}. $\bl$

The following result gives an estimate for $\|f(U)-f(V)\|$ for functions $f$ in the Zygmund class $\L_1$.

\begin{thm}
\label{oLu}
There exists a constant $c>0$ such that for every function $f\in\L_1$ and for arbitrary unitary operators $U$ and $V$ on Hilbert space
the following inequality holds:
$$
\|f(U)-f(V)\|\le c\,\|f\|_{\L_1}\left(2+\log_2\frac1{\|U-V\|}\right)\|U-V\|.
$$
\end{thm}

\Pf Again, as in the proof of Theorem \ref{uH}, we assume that $f=f_+$ and $N$ is the nonnegative integer satisfying \rf{N}. Using the notation introduced in the proof of Theorem \ref{uH}, we obtain
\begin{align*}
\left\|\sum_{n\le N}\big(f_n(U)-f_n(V)\big)\right\|&\le\sum_{n\le N}\big\|f_n(U)-f_n(V)\big\|\\[.2cm]
&\le\const\sum_{n\le N}2^n\|U-V\|\cdot\|f_n\|_{L^\be}\\[.2cm]
&\le\const(1+N)\|f\|_{\L_\a}\|U-V\|\\[.2cm]
&\le\const\|f\|_{\L_\a}\left(2+\log_2\frac1{\|U-V\|}\right)\|U-V\|.
\end{align*}

On the other hand,
\begin{align*}
\left\|\sum_{n>N}\big(f_n(U)-f_n(V)\big)\right\|&\le\sum_{n>N}2\|f_n\|_{L^\be}
\le\const\sum_{n>N}2^{-n}\|f\|_{\L_\a}\\[.2cm]
&\le\const2^{-N}\|f\|_{\L_\a}\le\const\|f\|_{\L_\a}\|U-V\|.\quad\bl
\end{align*}

In a similar way we can obtain an estimate for differences of order $n$ and functions in $\L_n$ for an arbitrary positive integer $n$.

Let us obtain now an analog of Theorem \ref{QAKg} for unitary operators. Let  $U$ be a unitary operator and let $A$ be a bounded
self-adjoint operator on a Hilbert space $\h$. Suppose that $f\in\L_\a$. By Theorem \ref{hou}, for every $u,\,v\in\h$,  the function
$$
t\mapsto f^{u,v}_{A,K}(t)\df\big(f(e^{{\rm i}tA}U)u,v\big)
$$
on $\R$ belongs to the space $\L_\a(\R)$. Thus for every $g\in B_1^{-\a}(\R)$, we can define the operator $\cR^g_{U,A}:\L_\a\to{\mathscr B}(\h)$
such that 
$$
\big(\big(\cR^g_{U,A}f\big)u,v\big)=\langle f^{u,v}_{A,K},g\rangle
$$
(here we identify the dual space $(B_1^{-\a}(\R))^*$ with $\L_\a(\R)$, see \S\,2.1).

\begin{thm}
\label{RUA}
Let $\a>0$. Then there exists $c>0$ such that for arbitrary unitary operator
$U$ and  a boundary self-adjoint operator $A$, and for every $g\in B_1^{-\a}(\R)$,
$$
\big\|\cR_{U,A}^g\big\|\le c\|g\|_{B_1^{-a}(\R)}\|A\|^\a.
$$
\end{thm}

\Pf Clearly,
\begin{align*}
\big|\big((\cR^g_{U,A}f)u,v\big)\big|&\le\const\big\|f^{u,v}_{A,K}\big\|_{\L_\a(\R)}\|g\|_{B_1^{-\a}(\R)}\\[.2cm]
&\le\const\|u\|\cdot\|v\|\cdot\|f\|_{\L_\a}\|g\|_{B_1^{-\a}(\R)}\|A\|^\a.\quad\bl
\end{align*}

\

\section{\bf The case of contractions}
\setcounter{equation}{0}
\label{s6}

\

In this section we obtain analogs of the results of \S\,4 and \S\,5 for contractions. Recall that if $T$ is a contraction on Hilbert space, it follows from von Neumann's inequality that the polynomial functional calculus $\f\mapsto f(T)$ extends to the disk-algebra
$C_A$ and $\|f(T)\|\le\|f\|_{C_A}$, $f\in C_A$.

\begin{thm}
\label{cH}
Let $0<\a<1$. Then there is a constant $c>0$ such that for 
every \lb$f\in(\L_\a)_+$ and for arbitrary contractions $T$ and $R$ on Hilbert space the following inequality holds:
$$
\|f(T)-f(R)\|\le c\,\|f\|_{\L_\a}\cdot\|T-R\|^\a.
$$
\end{thm}

\Pf
The proof of Theorem \ref{cH} is almost the same as the poof of Theorem \ref{uH}.  For $f\in(\L_\a)_+$, we use expansion
\rf{dia} and choose $N$ such that 
$$
2^{-N}<\|T-R\|\le2^{-N+1}.
$$
Then as in the proof of Theorem \ref{uH}, for $n\le N$, we estimate $\|f_n(T)-f_n(R)\|$ in terms of $\const 2^{-n}\|T-R\|$ (see \rf{BSc} and \rf{Bp}), while
for $n>N$ we use von Neumann's inequality to estimate $\|f_n(T)-f_n(R)\|$ in terms of $2\|f_n\|_{L^\be}$. The rest of the proof is the same. $\bl$

\begin{cor} 
\label{SK}
Let $f$ be a function in the disk algebra and $0<\a<1$.
Then the following two statements are equivalent:

{\em(i)} $\|f(T)-f(R)\|\le \const\|T-R\|^\a$ for all contractions $T$ and $R$,

{\em(ii)} $\|f(U)-f(V)\|\le \const\|U-V\|^\a$ for all unitary operators $U$ and $V$.
\end{cor}

{\bf Remark.} This corollary is also true for $\a=1$. This was proved
by Kissin and Shulman \cite{KS}.

\medskip






The following result is an analog of Theorem \ref{hou} for contractions.

\begin{thm}
\label{conh}
Let $m$ be a positive integer and $0<\a<m$. Then there exists a constant $c>0$ such that for every $f\in(\L_\a)_+$ and
for arbitrary contractions $T$ and $R$ on Hilbert space the following inequality holds:
$$
\left\|\sum_{k=0}^m(-1)^k\left(\begin{matrix}m\\k\end{matrix}\right)f\left(T+\frac{k}{m}(T-R)\right)
\right\|\le c\,\|f\|_{\L_\a}\|T-R\|^\a.
$$
\end{thm}

To prove Theorem \ref{conh}, we use the following analog of Lemma \ref{m}.

\begin{lem}
\label{mc}
Let $m$ be a positive integer and let $f$ be a function of class $\big(B^m_{\be1}\big)_+$. 
If $T$ and $R$ are contractions
on Hilbert space, then
\begin{align*}
&\sum_{k=0}^m(-1)^k\left(\begin{matrix}m\\k\end{matrix}\right)f\left(T+\frac{k}{m}(T-R)\right)\\[.2cm]
&=
\frac{m!}{m^m}\underbrace{\int\cdots\int}_{m+1}
(\dg^{m}f)(\z_1,\cdots,\z_{m+1})
\,d\E_1(\z_1)(T-R)\cdots(T-R)\,d\E_{m+1}(\z_{m+1}),
\end{align*}
where $\E_k$ is a semi-spectral measure of $T+\frac km(T-R)$.
\end{lem}

We conclude this section with an analog of Theorem \ref{oLu}.

\begin{thm}
\label{oLc}
There exists a constant $c>0$ such that for every function $f\in(\L_1)_+$ and for arbitrary contractions $T$ and $R$ on Hilbert space
the following inequality holds:
$$
\|f(T)-f(R)\|\le c\,\|f\|_{\L_1}\left(2+\log_2\frac1{\|T-R\|}\right)\|T-R\|.
$$ 
\end{thm}

\

\section{\bf Arbitrary moduli of continuity}
\setcounter{equation}{0}
\label{s7}

\

In this section we consider the problem of estimating $\|f(A)-f(B)\|$ for self-adjoint operators $A$ and $B$ and functions $f$ in the space $\L_\o$
(see \S\,2.2), where 
$\o$ is an arbitrary modulus of continuity. For simplicity, we give complete proofs for bounded self-adjoint operators. The case of unbounded self-adjoint operators will be considered in \cite{AP3}. We also obtain similar results for unitary operators and for contractions.

We have mentioned in the introduction that a Lipschitz function does not have to be operator Lipschitz and a continuously differentiable function does not have to be operator differentiable. On the other hand, we have proved in \S\,4 that a H\"older function of order $\a\in(0,1)$ must be operator H\"older of order $\a$ as well as a Zygmund function must be operator Zygmund. Moreover, the same is true for all classes $\L_\a$ with $\a>0$. This 
suggests an idea that the situation is similar with continuity properties of the Hilbert transform. In this section we consider the problem for which
moduli of continuity $\o$ the fact that $f\in\L_\o$ implies that $f$ belongs to the ``operator space $\L_\o$'', i.e.,
$$
\|f(A)-f(B)\|\le\const\o\big(\|A-B\|\big).
$$
We are going to compare this property with the fact that the Hilbert transform acts on $\L_\o$.

Given  a modulus of continuity $\o$, we define the function $\o_*$ by
$$
\o_*(x)=x\int_x^\be\frac{\o(t)}{t^2}\,dt,\quad x>0.
$$

\begin{thm}
\label{omsa}
There exists a constant $c>0$ such that for every
modulus of continuity $\o$, every $f\in \L_\o(\R)$ and for arbitrary self-adjoint operators $A$ and $B$,
the following inequality holds
$$
\|f(A)-f(B)\|\le c\,\|f\|_{\L_\o(\R)}\,\o_*\big(\|A-B\|\big).
$$
\end{thm}

\Pf Since $A$ and $B$ are bounded operators and their spectra are contained in $[a,b]$, we can replace a function $f\in\L_\o(\R)$ with the bounded function
$f_\flat$ defined by
\bay
\label{flat}
f_\flat(x)=\left\{\begin{array}{ll}f(b),&x>b,\\[.2cm]
f(x),&x\in[a,b],\\[.2cm]
f(a),&x<a.
\end{array}\right.
\ey
Clearly, $\|f_\flat\|_{\L_\o(\R)}\le\|f\|_{\L_\o(\R)}$. Thus we may assume that $f$ is bounded.

Let $N$ be an integer. We claim that 
\bay
\label{f}
f(A)-f(B)=\sum_{n=-\be}^N\big(f_n(A)-f_n(B)\big)+\big((f-f*V_N)(A)-(f-f*V_N)(B)\big),
\ey
and the series converges absolutely in the operator norm. Here $f_n=f*W_n+f*W_n^\sharp$ and the de la Vall\'ee Poussin type kernel $V_N$ is defined in \S\,2.1. Suppose that $M<N$. Indeed, it is easy to see that
\begin{align*}
f(A)-f(B)&-
\left(\sum_{n=M+1}^N\big(f_n(A)-f_n(B)\big)+\big((f-f*V_N)(A)-(f-f*V_N)(B)\big)\right)\\[.2cm]
&=\big((f-f*V_M)(A)-(f-f*V_M)(B)\big).
\end{align*}
Clearly, $f-f*V_M$ is an entire function of exponential type at most $2^{M+1}$. Thus it follows from \rf{Be} that
$$
\big\|(f-f*V_M)(A)-(f-f*V_M)(B)\big\|\le\const2^M\|f\|_{L^\be}\|A-B\|\to0\quad\mbox{as}\quad M\to-\be. 
$$

Suppose  now that $N$ is the integer satisfying \rf{AB}. It follows from Theorem \ref{Vn} that
\begin{align*}
\big\|(f-f*V_N)(A)&-(f-f*V_N)(B)\big\|\le2\|f-f*V_N\|_{L^\be}\\[.2cm]
&\le\const\|f\|_{\L_\o(\R)}\o\big(2^{-N}\big)\le\const\|f\|_{\L_\o(\R)}\o\big(\|A-B\|\big).
\end{align*}
On the other hand, it follows from Corollary \ref{Wnn} and from \rf{Be} that
\begin{align*}
\sum_{n=-\be}^N\|f_n(A)-f_n(B)\|&\le\const\sum_{n=-\be}^N2^n\|f_n\|_{L^\be}\|A-B\|\\[.2cm]
&\le\const\sum_{n=-\be}^N2^n\o\big(2^{-n}\big)\|f\|_{\L_\o(\R)}\|A-B\|\\[.2cm]
&=\const\sum_{k\ge0}2^{N-k}\o\big(2^{-N+k}\big)\|f\|_{\L_\o(\R)}\|A-B\|\\[.2cm]
&\le\const\left(\int_{2^{-N}}^\be\frac{\o(t)}{t^2}\,dt\right)\|f\|_{\L_\o(\R)}\|A-B\|\\[.2cm]
&=\const2^{N}\o_*\big(2^{-N}\big)\|f\|_{\L_\o(\R)}\|A-B\|\\[.2cm]
&\le\const\|f\|_{\L_\o(\R)}\o_*\big(\|A-B\|\big).
\end{align*}
The result follows now from the obvious inequality $\o(x)\le\o_*(x)$, $x>0$. $\bl$

\medskip

{\bf Remark.} Obviously, if  $\o_*(x)<\be$ for some $x>0$, then $\o_*(x)<\be$ for every $x>0$.
It follows easily from l'H\^opital's rule that in this case
$$
\lim_{x\to0}\o_*(x)=0.
$$
Moreover, in this case $\o_*$ is also a modulus of continuity. Indeed, it is easy to see that
$$
\o_*(x)=
\int_1^\infty\frac{\omega(sx)}{s^{2}}ds
$$
which implies that
$$
\o_*(x+y)\le\o_*(x)+\o_*(y),\quad x,~y\ge0
$$
and
$$
\o_*(x)\le\o_*(y),\quad 0\le x\le y.
$$

\medskip

Note that if the modulus of continuity $\o$ is bounded, then obviously, $\o_*(x)<\be$ for every $x>0$.
In the case when $A$ and $B$ are bounded self-adjoint operators and their spectra are contained in $[a,b]$,  we can 
replace $f$ with the function $f_\flat$ defined by \rf{flat} redefine the function $\o$ on $[b-a,\be)$ by putting $\o(x)=\o(b-a)$. 
Clearly, the modified modulus of continuity is bounded.

\begin{cor}
\label{mn}
Let $\o$ be a modulus of continuity such that 
$$
\o_*(x)\le\const\,\o(x),\quad x>0.
$$
Then for an arbitrary function $f\in\L_\o(\R)$ and for arbitrary self-adjoint operators $A$ and $B$ on Hilbert space the following inequality holds:
\bay
\label{n*}
\|f(A)-f(B)\|\le\const\|f\|_{\L_\o(\R)}\,\o\big(\|A-B\|\big).
\ey
\end{cor}

In the next result we do not pretend for maximal generality.

\begin{cor}
\label{ob}
Let $\o$ be a modulus of continuity such that $\o(2x)\le\vk\o(x)$ for some $\vk<2$ and all $x>0$.
Then $\o_*(x)\le\const\o(x)$ and
$$
\|f(A)-f(B)\|\le\const\|f\|_{\L_\o(\R)}\,\o\big(\|A-B\|\big)
$$
for arbitrary self-adjoint operators $A$ and $B$.
\end{cor}

\Pf It is easy to see that
$$
\o(t)\le\vk\left(\frac{t}{x}\right)^{\log_2\vk}\o(x),
$$
whenever $0<x\le t$. Thus
$$
\o_*(x)=x\int_x^\be\frac{\o(t)}{t^2}\,dt\le\vk x^{1-\log_2\vk}\o(x)\int_x^\be t^{\log_2\vk-2}\,dt\le
\frac{\vk}{1-\log_2\vk}\o(x).\quad\bl
$$

\medskip

{\bf Remark.} It is well known (see \cite{Z}, Ch. 3, Theorem 13.30) that if $\o$ is a modulus of continuity, then the Hilbert transform 
maps $\L_\o$ into itself if and only if
$$
\int_0^x\frac{\o(t)}{t}\,dt
+x\int_x^\be\frac{\o(t)}{t^2}\,dt\le\const\,\o(x),\quad x>0.
$$
It follows from Corollary \ref{mn} that if the Hilbert transform maps $\L_\o$ into itself, then \rf{n*} holds. However, the converse is false. For example, we can take a bounded modulus of continuity $\o$ such that $\o(x)$ is equivalent to 
$|\log x|^{-\a}$ near the origin and $\a>0$.

\medskip

In \cite{FN} it was proved that if $A$ and $B$ are self-adjoint operators on Hilbert space whose spectra
are contained in $[a,b]$ and $f$ is a continuous function on $[a,b]$, then
$$
\|f(A)-f(B)\|\le4\left(\log\left(\frac{b-a}{\|A-B\|}+1\right)+1\right)^2\o_f\big(\|A-B\|\big),
$$
where
$$
\o_f(\d)=\sup\big\{|f(x)-f(y)|:~x,\,y\in[a,b],~|x-y|<\d\big\}.
$$
The following corollary improves the result of Farforovskaya and Nikolskaya.

\begin{cor}
\label{FaN}
Suppose that $A$ and $B$ be self-adjoint operators with spectra in an interval $[a,b]$. Then for a continuous function $f$ on $[a,b]$ the following inequality holds:
$$
\|f(A)-f(B)\|\le\const\,\log\left(\frac{b-a}{\|A-B\|}+1\right)
\,\o_f\big(\|A-B\|\big).
$$
\end{cor}

\Pf Put $\o=\o_f$. Clearly, we may assume that $\o(x)=\o(b-a)$ for $x>a$.
Using the obvious inequality 
$$
\frac{\o(t)}t\le2\frac{\o(x)}x,\quad x\le t,
$$
we obtain
\begin{align*}
\o_*(x)&=x\int_x^\be\frac{\o(t)}{t^2}\,dt=x\int_x^{b-a}\frac{\o(t)}{t^2}\,dt+x\int_{b-a}^\be\frac{\o(t)}{t^2}\,dt\\[.2cm]
&\le2\o(x)\int_x^{b-a}\frac{dt}t+x\frac{\o(b-a)}{b-a}\le2\o(x)\log\frac{b-a}x+2\o(x)\\[.2cm]
&=2\o(x)\log\left(\frac{b-a}{x}+1\right).
\end{align*}
The result follows now from Theorem \ref{omsa}. $\bl$

\begin{cor}
\label{Lip} Let $f$ be a Lipschitz function on $\R$. Then for self-adjoint operators $A$ and $B$ with spectra in an interval $[a,b]$, the following
inequality holds
\bay
\label{uF}
\|f(A)-f(B)\|\le\const\|f\|_{\rm Lip}\,\log\left(\frac{b-a}{\|A-B\|}+1\right)
\|A-B\|.
\ey
\end{cor}

Note that a similar estimate can be obtained for bounded functions $f$ in the Zygmund class $\L_1(\R)$. This will be done at the end of the next section.

Inequality \rf{uF} improves the estimate 
$$
\|f(A)-f(B)\|\le\const\|f\|_{\rm Lip}\left(\log\left(\frac{b-a}{\|A-B\|}+1\right)+1\right)^2\|A-B\|.
$$
obtained  in \cite{F1} (see also \cite{F2}).

To conclude this section, we state analogs of Theorem \ref{omsa} for unitary operators and for contractions.

\begin{thm}
\label{omu}
There exists a constant $c>0$ such that for every
modulus of continuity $\o$, for every $f\in \L_\o$, and for arbitrary unitary operators $U$ and $V$,
the following inequality holds
$$
\|f(U)-f(V)\|\le c\,\|f\|_{\L_\o}\,\o_*(\|U-V\|).
$$
\end{thm}

\begin{thm}
\label{omc}
There exists a constant $c>0$ such that for every
modulus of continuity $\o$, for every $f\in\big( \L_\o\big)_+$, and for arbitrary contractions $T$ and $R$,
the following inequality holds
$$
\|f(T)-f(R)\|\le c\,\|f\|_{\L_\o}\,\o_*(\|T-R\|).
$$
\end{thm}

The proofs of Theorems \ref{omu} and \ref{omc} are similar to the proof of Theorem \ref{omsa}. Actually, they are even simpler, since we 
do not have to deal with convolutions with $W_n$ and $W_n^\sharp$ with negative $n$ which makes analogs of formula \rf{f} trivial.

\

\section{\bf Operator continuous functions and operator moduli of continuity}
\setcounter{equation}{0}
\label{s8}

\

In this section we introduce notions of operator continuous functions 
and uniformly operator continuous functions. We also define for a given continuous function on $\R$ the operator modulus of continuity associated with the function. We prove that a function is operator continuous if and only if it is uniformly operator continuous.

\medskip

{\bf Definition 1}. For  a continuous function $f$ on $\R$, we consider the map
\bay
\label{fA}
A\mapsto f(A)
\ey
defined on the set of (not necessarily bounded) self-adjoint operators.  We say that $f$ is {\it operator continuous}
if the map \rf{fA} is continuous at every (bounded or unbounded) self-adjoint operator $A$.

\medskip

This means that if $A$ is a (not necessarily bounded) self-adjoint operator, then for an arbitrary $\e>0$ there exists $\d>0$ such that
$\|f(A+K)-f(A)\|<\e$, whenever $K$ is a self-adjoint operator whose norm is less than $\d$. 

Note that it is easy to see that if $f$ is a continuous function on $\R$, then the map \rf{fA} is continuous at every bounded self-adjoint operator $A$.
Indeed, this is obvious for polynomials $f$. The result for arbitrary continuous functions follows from the Weirstrass theorem.

\medskip

{\bf Definition 2.} Let $f$ be a Borel function on $\R$. It is called {\it uniformly operator continuous} if for every $\e>0$ there exists $\d>0$ such that 
$\|f(A)-f(B)\|<\e$, whenever $A$ and $B$ are bounded self-adjoint operators such that $\|A-B\|<\d$.

\begin{thm}
\label{uoc}
Let $f$ be a bounded uniformly continuous function on $\R$. Then $f$ is uniformly operator continuous.
\end{thm}

\Pf Let $\o=\o_f$. Then $\o$ is a bounded modulus of continuity, and so $\o_*(x)<\be$, $x>0$. The result follows now from Theorem \ref{omsa} and the Remark
following that theorem. $\bl$

\medskip

{\bf Definition 3.}
Let $f$ be a continuous function on $\R$. Put
$$
\O_f(\d)\df\sup\big\|f(A)-f(B)\big\|,\quad\d>0,
$$
where the supremum is taken over all bounded self-adjoint operators $A$ and $B$ such that $\|A-B\|\le\d$. We say that $\O_f$ is the
{\it operator modulus of continuity of} $f$. 

\medskip

Note that it suffices to consider only operators $A$ and $B$ that are unitary equivalent to each other. Indeed, if $A$ and $B$ are self-adjoint operators
on a Hilbert space $\h$, we can define on the space $\h\oplus\h$ the self-adjoint operators $\A$ and $\B$ by
$$
\A=\left(\begin{matrix}A&\0\\\0&B\end{matrix}\right)\quad\mbox{and}\quad\B=\left(\begin{matrix}B&\0\\\0&A\end{matrix}\right).
$$
Obviously,
$$
\|\A-\B\|=\|A-B\|\quad\mbox{and}\quad\|f(\A)-f(\B)\|=\|f(A)-f(B)\|.
$$

We have by Theorem \ref{omsa},
$$
\o_f(\d)\le\O_f(\d)\le\const\o_*(\d),\quad\d>0.
$$

\medskip

\begin{thm}
\label{omc0}
Let $f$ be an operator continuous function. Then
$$
\lim_{\d\to0}\O_f(\d)=0,
$$
and so $f$ is uniformly operator continuous.
\end{thm}

\Pf Suppose that 
$$
\lim_{\d\to0}\O_f(\d)>\s>0.
$$
Then there are sequences of self-adjoint operators $\{A_j\}_{j\ge0}$ and $\{K_j\}_{j\ge0}$ on Hilbert space $\h$
such that $\|K_j\|<1/j$ and $\|f(A_j+K_j)-f(A_j)\|>\s$.
We define the operators $A$ and $R_n$ on $\ell^2(\h)$ by 
$$
A\left(\begin{matrix}h_0\\h_1\\h_2\\\vdots\end{matrix}\right)=\left(\begin{matrix}A_0h_0\\A_1h_1\\A_2h_2\\\vdots\end{matrix}\right)
\quad\mbox{and}\quad
R_n\left(\begin{matrix}h_0\\h_1\\h_2\\\vdots\end{matrix}\right)=\left(\begin{matrix}\0\\\vdots\\\0\\K_nh_n\\K_{n+1}h_{n+1}\\\vdots\end{matrix}\right).
$$
Clearly, $\|R_n\|\to0$ as $n\to0$, while $\|f(A+R_n)-f(A)\|>\s$ for $n\ge0$, and so the map \rf{fA} is not continuous at $A$. $\bl$

\medskip

{\bf Example.} Consider the function $g$ defined by $g(t)=|t|$, $t\in\R$. It was proved in \cite{Ka} that the function $g$ is not operator Lipschitz. 
It was observed in \cite{FN} that the function $g$ is not operator continuous. Let us show that 
$$
\O_g(\d)=\be\quad\mbox{for every}\quad\d>0,
$$
which will also imply that $g$ is not operator continuous. Indeed, suppose that $\O_g(\d_0)<\be$ for some $\d_0>0$. Since $g$ is homogeneous, it follows that 
$\O_g(\d)=\d\d_0^{-1}\O_g(\d_0)=\const\d$. However, this implies that $g$ is an operator Lipschitz function which contradicts the result of \cite{Ka}.


\begin{thm}
\label{Of}
Let  $A$ and $B$ be a pair of (not necessarily bounded) self-adjoint operators such that $A-B$ is bounded. Then
$$
\|f(A)-f(B)\|\le\O_f\big(\|A-B\|\big)
$$
for every continuous function $f$ on $\R$.
\end{thm}

To prove Theorem \ref{Of}, we need a couple of lemmata.

\begin{lem}
\label{sot} 
Let $f$ be a bounded continuous function on $\R$. Suppose that $A$ is a self-adjoint operator (not necessarily bounded) and $\{A_j\}_{j\ge0}$ is a sequence of bounded self-adjoint operators such that 
\bay
\label{Aj}
\lim_{j\to\be}\|A_ju-Au\|=0\quad\mbox{for every}\quad u\quad\mbox{in the domain of}\quad A.
\ey
Then
\bay
\label{sil}
\lim_{j\to\be}f(A_j)=f(A)\quad\mbox{in the strong operator topology}.
\ey
\end{lem}

\Pf We consider first the special case when $f(t)=(\l-t)^{-1}$, $\l\in\C\setminus\R$. Let $u$ be a vector in $\cd_A$, where $\cd_A$ 
denotes the domain of $A$.
Put $u_\l\df(\l I-A)^{-1}u$. Clearly, $u_\l\in \cd_A$ and
\begin{align*}
(\l I-A_j)^{-1}u&=(\l I-A_j)^{-1}(\l I -A)u_\l\\[.2cm]
&=u_\l+(\l I-A_j)^{-1}(A_ju-Au)\to u_\l\quad\mbox{as}\quad j\to\be.
\end{align*}
Since the linear combinations of such rational fractions are dense in the space $C_0(\R)$ of continuous functions on $\R$ vanishing at infinity,
it follows that \rf{sil} holds for an arbitrary function $f$ in $C_0(\R)$.

Suppose now that $f$ is an arbitrary bounded continuous function on $\R$. By subtracting from $f$ a continuous function with compact support, we may assume that $f$ vanishes on $[-1,1]$. Then there exists a function $g$ in $C_0(\R)$ such that $f(t)=tg(t)$, $t\in\R$. Let $u\in \cd_A$. We have
\begin{align}
\label{str}
\nonumber
f(A_j)u&=g(A_j)A_ju=g(A_j)Au+g(A_j)(A_ju-Au)\\[.2cm]
&\to g(A)Au=f(A)u\quad\mbox{as}\quad j\to\be.\quad\bl
\end{align}

\begin{lem}
\label{CR}
Let $f$ be a continuous function on $\R$ such that $|f(t)|\le\const(1+|t|)$, $t\in\R$ and let $A$ and $\{A_j\}_{j\ge0}$ be as in Lemma
{\em\ref{sot}}. Then
$$
\lim_{j\to\be}\|f(A_j)u-f(A)u\|=0\quad\mbox{for every}\quad u\in \cd_A.
$$
\end{lem}

\Pf As in the proof of Lemma \ref{sil}, we may assume that $f$ vanishes on $[-1,1]$ and define the continuous function $g$ by $f(t)=tg(t)$, $t\in\R$.
It follows now from Lemma \ref{sil} that \rf{str} holds for every $u\in \cd_A$. $\bl$

\medskip

{\bf Proof of Theorem \ref{Of}.} Clearly, if $\O_f(\d)<\be$ for some $\d>0$, it follows that $f$ satisfies the hypotheses of Lemma \ref{CR}. Let $K=B-A$. Then $K$ is a bounded self-adjoint operator. Put
$$
A_j\df E_A\big([-j,j]\big)A.
$$
Clearly, \rf{Aj} holds. It follows easily from Lemma \ref{str} that
$$
\|f(A+K)-f(A)\|\le\limsup_{j\to\be}\|f(A_j+K)-f(A_j)\|\le\O_f\big(\|K\|\big).\quad\bl
$$


\begin{cor}
\label{sled}
Let $f$ be continuous function on $\R$. Then $f$ is operator continuous if and only if it is uniformly operator continuous.
\end{cor}

We conclude this section with an estimate for the operator modulus of continuity of a bounded function in the Zygmund class $\L_1(\R)$. The proof
of the following theorem is similar to the proof of Theorem 3.4 of Ch. 2 of \cite{Z}.

\begin{thm}
\label{Zyg}
Let $f$ be a bounded function in $\L_1(\R)$. Then there exists $c>0$ such that
$$
\O_f(\d)\le c\,\d\log\frac2\d\quad\mbox{for}\quad\d\le1.
$$
\end{thm}

\Pf By Corollary \ref{Z}, there is a constant $c_1$ such that
$$
\big\|f(A+2K)-2f(A+K)+f(K)\big\|\le c_1\|f\|_{\L_1(\R)}\|K\|.
$$
It is easy to see that
\begin{align*}
\big\|f(A+K)-f(A)\big\|&\le\frac12\big\|f(A+2K)-2f(A+K)+f(K)\big\|\\[.2cm]
&+\frac12\big\|f(A+2K)-f(A)\big\|.
\end{align*}
It follows that
$$
\O_f(t/2)\le\frac{c_1}4\|f\|_{\L_1(\R)}t+\frac12\O_f(t),
$$
and so
$$
2^{k-1}\O_f\big(2^{-k}t\big)-2^{k-2}\O_f\big(2^{1-k}t\big)\le\frac{c_1}4\|f\|_{\L_1(\R)}t, \quad\mbox{whenever}\quad k\ge1.
$$
Substituting $t=t_0\df\frac4{c_1}\|f\|^{-1}_{\L_1(\R)}\|f\|_{L^\be}$, and keeping in mind the trivial estimate $\O_f(t)\le2\|f\|_{L^\be}$, $t>0$,
we obtain
$$
2^{n-1}\O_f\big(2^{-n}t_0\big)\le(n+1)\|f\|_{L^\be}.
$$
Hence, for $t=2^{-n}t_0$, $n\ge0$, we have
$$
\O_f(t)\le\frac{c_1}2\|f\|_{\L_1(\R)}t\log_2\left(\frac{8\|f\|_{L^\be}}{c_1\|f\|_{\L_1(\R)}t}\right)
$$
Therefore
$$
\O_f(t)\le c_1\|f\|_{\L_1(\R)}t\log_2\left(\frac{8\|f\|_{L^\be}}{c_1\|f\|_{\L_1(\R)}t}\right)\quad\mbox{for}\quad t\le\frac{t_0}2
$$
and $\O_f(t)\le2\|f\|_{L^\be}$ for $t\ge t_0/2$. $\bl$

\

\section{\bf A universal family of self-adjoint operators}
\setcounter{equation}{0}
\label{s9}

\

In this section we construct a universal family of (unbounded) self-adjoint operators $\{A_t\}_{t\ge0}$ such that
the operators $A_t$ have purely point spectra and
$$
\O_f(t)=\|f(A_t)-f(A_0)\|,\quad t>0,
$$
for every continuous function $f$. In particular, $\|A_t-A_0\|=t$, $t\ge0$.
Moreover, the operators $A_t$, $t\ge0$, are unitarily equivalent to each other.

Denote by ${\frak K}$ the set of finite rank self-adjoint operators on Hilbert space and let
${\frak K}_0$ be a countable dense subset of ${\frak K}$.

\begin{lem}
\label{siln}
Suppose that $\{A_j\}$ be a sequence of bounded self-adjoint operators that converges to $A$ in the strong operator topology.
Then $f(A_j)\to f(A)$ strongly for an arbitrary continuous function $f$.
\end{lem}

\Pf The conclusion of the lemma is trivial if $f$ is a polynomial. It remains to approximate $f$ by polynomials uniformly on 
$\big[-\sup_j\|A_j\|,\sup_j\|A_j\|\big]$. $\bl$

\begin{cor}
\label{sle}
Let $f\in C(\R)$ and $t>0$. Then
$$
\O_f(t)=\sup\big\{\|B-A\|:~A,\,B\in{\frak K}_0(\h),~\|B-A\|<t\big\}.
$$
\end{cor}

\Pf Clearly, we have to verify that the left-hand side is less than or equal to the right-hand side. Let $A$ and $B$ be bounded self-adjoint operators such that
$\|A-B\|<t$. Let $\{A_j\}$ and $\{K_j\}$ be sequences of operators in ${\frak K}_0$ such that $A_j\to A$, $K_j\to B-A$ in the strong operator topology,
and $\|K_j\|\le\|B-A\|$ for all $j$. By Lemma \ref{siln}, $f(A_j)\to f(A)$ and $f(A_j+K_j)\to f(B)$ strongly. Hence,
$$
\|f(B)-f(A)\|\le\liminf_{j\to\be}\|f(A_j+K_j)-f(A_j)\|
$$
which implies the desired inequality. $\bl$

Suppose that $\{R_j\}_{j=1}^\be$ is an enumeration of ${\frak K}_0$. For given $j\ge1$ and $t>0$ we consider the set
$$
{\frak K}_{jt}\df\big\{A\in{\frak K}_0:~\|A-R_j\|<t\big\}
$$
and let $\big\{R_{jk}^{(t)}\big\}_{k=1}^\be$ be an enumeration of ${\frak K}_{jt}$. Put $R_{jt}^{(0)}\df R_j$.

We can define now a universal family $\{A_t\}_{t\ge0}$ by
\bay
\label{os}
A_t\df\bigoplus_{j=1}^\be\bigoplus_{k=1}^\be R_{jk}^{(t)}.
\ey

\begin{thm}
\label{ufa}
The operators $A_t$ are pairwise unitarily equivalent. Each operator $A_t$ has purely point spectrum.
Moreover, for every continuous function $f$ on $\R$, we have
$$
\|f(A_t)-f(A_0)\|=\O_f(t),\quad t>0.
$$
\end{thm}

\Pf It is easy to see that each operator in ${\frak K}_0$ occurs in the orthogonal sum on the right of \rf{os} infinitely many times and 
each operator in the orthogonal sum on the right of \rf{os} belongs to ${\frak K}_0$. Thus $A_t$ is unitarily equivalent to $A_0$ for all $t>0$.

We have
$$
\|f(A_t)-f(A_0)\|=\sup_{j,k}\left\|f\left(R_{jk}^{(t)}\right)-f\left(R_{jk}^{(0)}\right)\right\|=\O_f(t)
$$
by Corollary \ref{sle}. $\bl$



\

\section{\bf Commutators and quasicommutators}
\setcounter{equation}{0}
\label{s10}

\

In this section we obtain estimates for the norm of {\it quasicommutators} $f(A)Q-Qf(B)$ in terms of $\|AQ-QB\|$ for self-adjoint operators $A$ and $B$
and a bounded operator $Q$. We assume for simplicity that $A$ and $B$ are bounded. However, we obtain estimates that do not depend
on the norms of $A$ and $B$. In \cite{AP3} we will  consider the case of not necessarily bounded operators $A$ and $B$.
Note that in the special case $A=B$ this problem turns into the problem of estimating the norm of commutators $f(A)Q-Qf(A)$ in terms of 
$\|AQ-QA\|$. On the other hand, in the special case $Q=I$ the problem turns into the problem of estimating $\|f(A)-f(B)\|$ in terms
$\|A-B\|$. 

Note that similar results can be obtained for unitary operators and for contractions.

Birman and Solomyak (see \cite{BS5}) discovered the following formula
$$
f(A)Q-Qf(B)=\iint\frac{f(x)-f(y)}{x-y}\,dE_A(x)(AQ-QB)\,dE_B(y),
$$
whenever $f$ is a function such that the divided difference ${\frak D}f$ is a Schur multiplier with respect to the spectral measures $E_A$ and $E_B$.

We could use this formula to obtain estimates of quasicommutators as we have done in the case of functions of perturbed operators.
However, we are going to reduce estimates of quasicommutators to those of functions of perturbed operators. For this purpose we obtain estimates that compare different moduli of continuities (the operator modulus of continuity, 
the (quasi)commutator modulus of continuity, etc).

We start with the case of operator Lipschitz functions.

The following theorem compares different operator Lpschitz norms
and (quasi)com\-mutator Lipschitz norms. The fact that they are equivalent is well-known, see \cite{KS}. The following theorem says that all those norms are equal.

\begin{thm}
\label{sr}
Let $f$ be a continuous function on $\R$. The following are equivalent:

{\em(i)} $\|f(A)-f(B)\|\le\|A-B\|$ for arbitrary self-adjoint operators $A$ and $B$;

{\em(ii)} $\|f(A)-f(B)\|\le\|A-B\|$ for all pairs of unitarily equivalent self-adjoint operators
$A$ and $B$;

{\em(iii)} $\|f(A)R-Rf(A)\|\le\|AR-RA\|$ for arbitrary self-adjoint operators $A$ and $R$;

{\em(iv)} $\|f(A)R-Rf(A)\|\le\|AR-RA\|$ for all self-adjoint operators $A$ and bounded operators $R$;

{\em(v)} $\|f(A)R-Rf(B)\|\le\|AR-RB\|$ for arbitrary self-adjoint operators $A$ and $B$ and an arbitrary bounded operator $R$.
\end{thm}

\Pf The implication (i)$\imp$(ii) is obvious.

Let us show that (ii)$\imp$(iii). Put $B=\exp(-{\rm i}tR)A\exp({\rm i}tR)$. Clearly, $B$ is unitarily equivalent to $A$ and 
$f(B)=\exp(-{\rm i}tR)f(A)\exp({\rm i}tR)$. Thus
$$
\big\|f(A)-\exp(-{\rm i}tR)f(A)\exp({\rm i}tR)\big\|\le\big\|A-\exp(-{\rm i}tR)A\exp({\rm i}tR)\big\|\quad\mbox{for all}\quad t\in\R.
$$
It remains to observe that
$$
\lim_{t\to0}\frac{\|f(A)-\exp(-{\rm i}tR)f(A)\exp({\rm i}tR)\|}{|t|}=\|f(A)R-Rf(A)\|
$$
and
$$
\lim_{t\to0}\frac{\|A-\exp(-{\rm i}tR)A\exp({\rm i}tR)\|}{|t|}=\|AR-RA\|.
$$

To prove that (iii)$\imp$(iv), we consider the following self-adjoint operators
$$
\A=\left(\begin{matrix}A&\0\\[.2cm]\0&A\end{matrix}\right)\quad\mbox{and}\quad
\cR=\left(\begin{matrix}\0&R\\[.2cm]R^*&\0\end{matrix}\right).
$$
It is easy to see that
$$
f(\A)\cR=\left(\begin{matrix}\0&f(A)R\\[.2cm]f(A)R^*&\0\end{matrix}\right)\quad\mbox{and}\quad
\cR f(\A)=\left(\begin{matrix}\0&R f(A)\\[.2cm]R^*f(A)&\0\end{matrix}\right).
$$
Hence,
$$
\|f(\A)\cR-\cR f(\A)\|=\max\big\{\|f(A)R-Rf(A)\|,~\|f(A)R^*-R^*f(A)\|\big\}
$$
and
$$
\|\A\cR-\cR \A\|=\max\big\{\|AR-RA\|,~\|AR^*-R^*A\|\big\}=\|AR-RA\|.
$$
It follows that
$$
\|f(A)R-Rf(A)\|\le\|f(\A)\cR-\cR f(\A)\|\le\|\A\cR-\cR \A\|=\|AR-RA\|.
$$

The implication (v)$\imp$(i) is trivial; it suffices to put $R=I$. 

To complete the proof, it remains to show that (iv)$\imp$(v). Let us first consider the special case when $A$ and $B$
are unitarily equivalent, i.e., $A=U^*BU$ for a unitary operator $U$ and we prove that 
$$
\|U^*f(B)UR-Rf(B)\|\le\|U^*BUR-RB\|.
$$
This is equivalent to the inequality
$$
\|f(B)UR-URf(B)\|\le\|BUR-URB\|
$$
which holds by (iv).

Now we consider the case of arbitrary self-adjoint operators $A$ and $B$. Put
$$
\A=\left(\begin{matrix}A&\0\\[.2cm]\0&B\end{matrix}\right),\quad
\B=\left(\begin{matrix}B&\0\\[.2cm]\0&A\end{matrix}\right),\quad\mbox{and}\quad
\cR=\left(\begin{matrix}R&\0\\[.2cm]\0&R^*\end{matrix}\right).
$$
Then $\A$ and $\B$ are unitarily equivalent. We have
$$
f(\A)\cR=\left(\begin{matrix}f(A)R&\0\\[.2cm]\0&f(B)R^*\end{matrix}\right)\quad\mbox{and}\quad
\cR f(\B)=\left(\begin{matrix}Rf(B)&\0\\[.2cm]\0&R^*f(A)\end{matrix}\right).
$$
Hence,
$$
\|f(\A)\cR-\cR f(\B)\|=\max\big\{\|f(A)R-Rf(B)\|,~\|f(B)R^*-R^*f(A)\|\big\}
$$
and
$$
\|\A\cR-\cR\B\|=\max\big\{\|AR-RB\|,~\|BR^*-R^*A\|\big\}=\|AR-RB\|.
$$
It follows that
$$
\|f(A)R-Rf(B)\|\le\|f(\A)\cR-\cR f(\B)\|\le\|\A\cR-\cR\B\|=\|AR-RB\|.\quad\bl
$$

In \S\,8 for a continuous function $f$ on $\R$ we have defined the operator modulus of continuity $\O_f$. 
We define here 3 other version of moduli of continuity in terms of commutators and quasicommutators.

Let $f$ be a continuous function on $\R$. For $\d>0$, put
\begin{align*}
\O_f^{[1]}(\d)&\df\sup\big\{\|f(A)R-Rf(A)\|:~A~\mbox{and}~R~\mbox{are self-adjoint},~\|R\|=1\big\};\\[.2cm]
\O_f^{[2]}(\d)&\df\sup\big\{\|f(A)R-Rf(A)\|:~A~\mbox{is self-adjoint},~\|R\|=1\big\};\\[.2cm]
\O_f^{[3]}(\d)&\df\sup\big\{\|f(A)R-Rf(B)\|:~A~\mbox{and}~B~\mbox{are self-adjoint},~\|R\|=1\big\}.
\end{align*}
Obviously, $\O_f^{[1]}\le\O_f^{[2]}\le\O_f^{[3]}$ and $\O_f\le\O_f^{[3]}$.

\begin{thm}
\label{mcc}
Let $f$ be a continuous function on $\R$. Then
$$
\O_f\le\O_f^{[1]}=\O_f^{[2]}=\O_f^{[3]}\le2\O_f.
$$
\end{thm}

\Pf The inequality $\O_f^{[2]}\le\O_f^{[1]}$ can be proved in the same way as the implication (iii)$\imp$(iv) in the proof of Theorem 
\ref{sr}. The inequality $\O_f^{[3]}\le\O_f^{[2]}$ can be proved in the same way as the implication (iv)$\imp$(v) in the proof of Theorem \ref{sr}. It remains to prove that $\O_f^{[1]}\le2\O_f$. We need two lemmata.

\begin{lem}
\label{pl}
Let $X$ and $Y$ be bounded operators. Then
$$
\|XY^n-Y^nX\|\le n\|Y\|^{n-1}\|XY-YX\|.
$$
\end{lem}

\Pf We have
$$
\|XY^n-Y^nX\|\le\left\|\sum_{k=1}^nY^{k-1}(XY-YX)Y^{n-k}\right\|\le n\|Y\|^{n-1}\|XY-YX\|.\quad\bl
$$

\begin{lem}
\label{vl}
Let $T$  be a self-adjoint operator such that $\|T\|<1$ and let $X$ be a bounded operator. Then
$$
\big\|(I-T^2)^{1/2}X-X(I-T^2)^{1/2}\big\|\le\frac{\|T\|\cdot\|XT-TX\|}{(1-\|T\|^2)^{1/2}}.
$$
\end{lem}

\Pf Let $a_n\df(-1)^{n-1}\left(\begin{matrix}1/2\\n\end{matrix}\right)$. Then $a_n>0$ and
$(1-t^2)^{1/2}=1-\sum\limits_{n=1}^\be a_nt^{2n}$. Thus
\begin{align*}
\big\|(I-T^2)^{1/2}X-X(I-T^2)^{1/2}\big\|&=\left\|\sum_{n=1}^\be a_n \big(XT^{2n}-T^{2n}X\big)\right\|\\[.2cm]
&\le\|XT-TX\|\sum_{n=1}^\be2na_n\|T\|^{2n-1}=\frac{\|T\|\cdot\|XT-TX\|}{(1-\|T\|^2)^{1/2}}
\end{align*}
by Lemma \ref{pl}. $\bl$

Let us complete the proof of Theorem \ref{mcc}. Let $R$ be a self-adjoint contraction and $\t\in(0,1)$.
Consider the operators
$$
\A=\left(\begin{matrix}A&\0\\[.2cm]\0&A\end{matrix}\right)\quad\mbox{and}\quad
\cU=\left(\begin{matrix}\t R&(I-\t^2R)^{1/2}\\[.2cm]-(I-\t^2R)^{1/2}&\t R\end{matrix}\right).
$$
Clearly, $\cU$ is a unitary operator. We have
$$
f(\A)\cU=\left(\begin{matrix}\t f(A)R&f(A)(I-\t^2R^2)^{1/2}\\[.2cm]-f(A)(I-\t^2R^2)^{1/2}&\t f(A)R\end{matrix}\right)
$$
and
$$
\cU f(\A)=\left(\begin{matrix}\t Rf(A)&(I-\t^2R^2)^{1/2}f(A)\\[.2cm]-(I-\t^2R^2)^{1/2}f(A)&\t Rf(A)\end{matrix}\right).
$$
Clearly, 
$$
\|f(\A)\cU-\cU f(\A)\|\ge\t\|f(A)R-Rf(\A)\|
$$
and
\begin{align*}
\|\A\cU-\cU\A\|&\le\t\|AR-RA\|+\big\|A\big(I-\t^2R^2\big)^{1/2}-\big(I-\t^2R^2\big)^{1/2}A\big\|\\[.2cm]
&\le\big(\t+\t^2(1-\t^2)^{-1/2}\big)\|AR-RA\|
\end{align*}
by Lemma \ref{vl}. Hence,
\begin{align*}
\|f(A)R-Rf(\A)\|&\le\t^{-1}\|f(\A)\cU-\cU f(\A)\|=\t^{-1}\|\cU^*f(\A)\cU-f(\A)\|\\[.2cm]
&\le\t^{-1}\O_f\big(\big\|\cU^*\A\cU-\A\big\|\big)=\t^{-1}\O_f\big(\big\|\A\cU-\cU\A\big\|\big)\\[.2cm]
&\le\t^{-1}\O_f\Big(\big(\t+\t^2(1-\t^2)^{-1/2}\big)\|AR-RA\|\Big).
\end{align*}
Taking $\t=1/2$, we obtain
\begin{align*}
\|f(A)R-Rf(\A)\|&\le2\O_f\left(\left(\frac12+\frac1{2\sqrt{3}}\right)\|AR-RA\|\right)\le2\O_f\big(\|AR-RA\|\big).\quad\bl
\end{align*}

\medskip

{\bf Remark.} It can be shown that there exist a uniformly continuous function $f$ and a positive number $\d$ such that
$\O_f(\d)<\O_f^{[1]}(\d)$. This will be shown in \cite{AP3}.

\medskip

Now we can deduce from Theorem  \ref{mcc} analogs of Theorems \ref{saH} and \ref{omsa} for quasicommutators.

\begin{thm}
\label{fcc}
Let $0<\a<1$. Then there exists $c>0$ such that for every $f\in\L_\a(\R)$, for arbitrary self-adjoint operators $A$ and $B$ and a bounded operator $R$ the following inequality holds:
$$
\|f(A)R-Rf(B)\|\le c\,\|f\|_{\L_\a(\R)}\|AR-RB\|^\a\|R\|^{1-\a}.
$$
\end{thm}

\Pf Clearly, we may assume that $Q\ne\0$.
By Theorems \ref{saH} and \ref{mcc},
\begin{align*}
\|f(A)R-Rf(B)\|&=\|R\|\cdot
\left\|f(A)\left(\frac1{\|R\|}R\right)-\left(\frac1{\|R\|}R\right)f(A)\right\|\\[.2cm]
&\le\const\|f\|_{\L_\a(\R)}\|R\|
\left\|\frac1{\|R\|}\big(f(A)R-Rf(A)\big)\right\|^\a\\[.2cm]
&=\const\|f\|_{\L_\a(\R)}\|AR-RB\|^\a\|R\|^{1-\a}.\quad\bl
\end{align*}

\begin{thm}
\label{oqc}
There exists $c>0$ such that for every
modulus of continuity $\o$, for every $f\in\L_\o(\R)$, for arbitrary self-adjoint operators $A$ and $B$, and a bounded nonzero
operator $R$ the following inequality holds:
$$
\|f(A)R-Rf(B)\|\le c\|R\|\,\,\o_*\!\left(\frac{\left\|\big(f(A)R-Rf(A)\big)\right\|}{\|R\|}\right).
$$
\end{thm}

The proof of Theorem \ref{oqc} is the same as the proof of Theorem \ref{fcc}.

\

\section{\bf Higher order moduli of continuity}
\setcounter{equation}{0}
\label{s11}

\



In this section we obtain norm estimates for finite differences
$$
\big(\D_K^mf\big)(A)\df\sum_{j=0}^m(-1)^{m-j}\left(\begin{matrix}m\\j\end{matrix}\right)f\big(A+jK\big)
$$
for functions $f\in\L_{\o,m}(\R)$ and self-adjoint operators $A$ and $K$. For simplicity, we give proofs in this paper in the case of bounded operators and bounded functions $f$. Note that our estimate will not depend on the $L^\be$ norm of $f$, nor on the operator norm of $A$.
In \cite{AP3} we consider the case of an arbitrary (not  necessarily bounded) self-adjoint operator $A$ (though $K$ still must be bounded) and an arbitrary function $f\in\L_{\o,m}(\R)$.

We also obtain similar results for unitary operators and for contractions.

Let
$\o$ be a nondecreasing function on $(0,\be)$ such that 
\bay
\label{om}
\lim_{x\to0}\o(x)=0\quad\mbox{and}\quad
\o(2x)\le2^m\o(x)\quad\mbox{for}\quad x>0.
\ey

Recall that $\L_{\o,m}(\R)$ is the space of continuous functions $f$ on $\R$
satisfying
$$
\|f\|_{\L_{\o,m}(\R)}\df\sup\limits_{t>0}\frac{\|\D^m_tf\|_{L^\infty}}{\o(t)}<+\infty.
$$

Given a nondecreasing function $\o$ satisfying \rf{om}, we define the function $\o_{*,m}$ by
$$
\o_{*,m}(x)=x^m\int_x^\be\frac{\o(t)}{t^{m+1}}\,dt=\int_1^\be\frac{\o(sx)}{s^{m+1}}\,dx.
$$

\begin{thm}
\label{oon}
Let $m$ be a positive integer.  Then there is a positive number $c$ such that for an arbitrary
nondecreasing function $\o$ on $(0,\be)$ satisfying {\em\rf{om}}, an arbitrary
bounded function $f$ in $\L_{\o,m}(\R)$, and arbitrary bounded self-adjoint operators $A$ and 
$K$ on Hilbert space the following inequality holds:
$$
\left\|\big(\D_K^mf\big)(A)\right\|\le c\,\|f\|_{\L_{\o,m}(\R)}\,\o_{*,m}\big(\|K\|\big).
$$
\end{thm}

\Pf As in the proof of Theorem \ref{omsa}, we can easily see that 
$$
\big(\D_K^mf\big)(A)=\sum_{n=-\be}^N\big(\D_K^mf_n\big)(A)+\big(\D_K^m(f-f*V_N)\big)(A),
$$
where as before, $f_n=f*W_n+f*W_n^\sharp$.

Suppose that $N$ is the integer satisfying \rf{K}. By Theorem \ref{mnn},
\begin{align*}
\big\|\big(\D_K^m(f-f*V_N)\big)(A)\big\|&\le\const\|f-f*V_N\|_{L^\be}\\[.2cm]
&\le\const\|f\|_{\L_{\o,m}(\R)}\o\big(2^{-N}\big)
\le\const\|f\|_{\L_{\o,m}(\R)}\o_{*,m}\big(||K\|\big).
\end{align*}
On the other hand, it follows from Lemma \ref{m}, \rf{Boke}, and Corollary \ref{Wnm} that
$$
\big\|\big(\D_K^mf_n\big)(A)\big\|\le\const2^{mn}\|f_n\|_{L^\be}\|K\|^m
\le\const\|f\|_{\L_{\o,m}(\R)}2^{mn}\o\big(2^{-n}\big)\|K\|^m.
$$
Thus
\begin{align*}
\sum_{n=-\be}^N\big\|\big(\D_K^mf_n\big)(A)\big\|&\le\const\sum_{n=-\be}^N\|f\|_{\L_{\o,m}(\R)}2^{mn}\o\big(2^{-n}\big)\|K\|^m\\[.2cm]
&=\sum_{k\ge0}2^{(N-k)m}\o\big(2^{N-k}\big)\|f\|_{\L_{\o,m}(\R)}\|K\|^m\\[.2cm]
&\le\const\left(\int_{2^{-N}}^\be\frac{\o(t)}{t^{m+1}}\,dt\right)\|f\|_{\L_{\o,m}(\R)}\|K\|^m\\[.2cm]
&=\const2^{-Nm}\o_{*,m}\big(2^{-N}\big)\|f\|_{\L_{\o,m}(\R)}\|K\|^m\\[.2cm]
&\le\const\|f\|_{\L_{\o,m}(\R)}\o_{*,m}\big(\|K\|\big).
\end{align*}
This completes the proof. $\bl$

\begin{cor}
\label{hvk}
Let $\o$ be a positive nondecreasing function on $(0,\be)$ such that $\lim\limits_{x\to0}\o(x)=0$ and
$\o(2x)\le\vk\o(x)$ for some $\vk<2^m$ and all $x>0$. Then for $x>0$, we have $\o_{*,m}(x)\le\const\o(x)$ and so
$$
\left\|\big(\D_K^mf\big)(A)\right\|\le c\,\|f\|_{\L_{\o,m}(\R)}\,\o\big(\|K\|\big).
$$
\end{cor}

The proof of Corollary \ref{hvk} is the similar to the proof of Corollary
\ref{ob}.

\begin{cor}
\label{M}
Suppose that under the hypotheses of Theorem {\em\ref{oon}} $\|f\|_{L^\be}\le M$.  Then for the function $\o_{m,M}$ defined by
$$
\o_{m,M}(x)=x^m\int_x^\be\frac{\min\big(2M,\|f\|_{\L_{\o,m}}(t)\big)}{t^{m+1}}\,dt,
$$
the following inequality holds:
$$
\big\|\big(\D_K^mf\big)(A)\big\|\le\const\|f\|_{\L_{\o,m}}\o_{m,M}\big(\|K\|\big).
$$
\end{cor}

The following analogs of Theorem \ref{oon} for unitary operators and for contractions can be proved in a similar way.

\begin{thm}
\label{uom}
Let $m$ be a positive integer. Then there exists a constant $c>0$ such that for every
nondecreasing function $\o$ on $(0,\be)$ satisfying {\em\rf{om}},  for every $f\in\L_{\o,m}$, and
for an arbitrary unitary operator $U$ and an arbitrary bounded self-adjoint operator $A$ on Hilbert space, the following inequality holds:
$$
\left\|\sum_{k=0}^n(-1)^k\left(\begin{matrix}m\\k\end{matrix}\right)f\big(e^{{\rm i}kA}U\big)
\right\|\le c\,\|f\|_{\L_{\o,m}}\o_{*,m}\big(\|A\|\big).
$$
\end{thm}

\begin{thm}
\label{co}
Let $m$ be a positive integer. Then there exists a constant $c>0$ such that for every
nondecreasing function $\o$  on $(0,\be)$ satisfying {\em\rf{om}}, for every \lb$f\in(\L_{\o,m})_+$, and
for arbitrary contractions $T$ and $R$ on Hilbert space the following inequality holds:
$$
\left\|\sum_{k=0}^m(-1)^k\left(\begin{matrix}m\\k\end{matrix}\right)f\left(T+(-1)^k\frac{k}{n}(T-R)\right)
\right\|\le c\,\|f\|_{\L_{\o,m}}\o_{*,m}\big(\|T-R\|\big).
$$
\end{thm}

\

\

\noindent
\begin{tabular}{p{9cm}p{15cm}}
A.B. Aleksandrov & V.V. Peller \\
St-Petersburg Branch & Department of Mathematics \\
Steklov Institute of Mathematics  & Michigan State University \\
Fontanka 27, 191023 St-Petersburg & East Lansing, Michigan 48824\\
Russia&USA
\end{tabular}

\end{document}